\documentclass[11pt,a4paper]{article}
\usepackage[margin=26mm]{geometry}
\usepackage[T1]{fontenc}
\usepackage[utf8]{inputenc}
\usepackage{lmodern,amsmath,amssymb,amsthm,mathtools,microtype}
\usepackage{booktabs,array,enumitem,xurl,longtable}
\usepackage{mathrsfs,needspace}
\usepackage{tikz}
\usetikzlibrary{arrows.meta,positioning}
\tikzset{
 dynkin edge/.style={line width=.8pt},
 dynkin plus/.style={circle,draw,fill=black,inner sep=0pt,minimum size=5.5pt},
 dynkin minus/.style={circle,draw,fill=white,inner sep=0pt,minimum size=5.5pt},
 dynkin neutral/.style={circle,draw,fill=gray!35,inner sep=0pt,minimum size=5.5pt}
}
\newcommand{\AiiDiagram}{%
\begin{tikzpicture}[x=.7cm,y=.7cm,baseline=-.6ex]
\node[dynkin plus] (a) at (0,0) {};
\node[dynkin minus] (b) at (1,0) {};
\draw[dynkin edge] (a)--(b);
\end{tikzpicture}}
\newcommand{\ArDiagram}{%
\begin{tikzpicture}[x=.7cm,y=.7cm,baseline=-.6ex]
\node[dynkin plus] (a) at (0,0) {};
\node[dynkin minus] (b) at (1,0) {};
\node[fill=white,inner sep=2pt] (c) at (2,0) {$\cdots$};
\node[dynkin neutral] (d) at (3,0) {};
\draw[dynkin edge] (a)--(b)--(c)--(d);
\end{tikzpicture}}
\newcommand{\DtreeDiagram}{%
\begin{tikzpicture}[x=.7cm,y=.7cm,baseline=-.6ex]
\node[dynkin minus] (u) at (-1,.55) {};
\node[dynkin minus] (v) at (-1,-.55) {};
\node[dynkin plus] (c) at (0,0) {};
\node[dynkin minus] (a) at (1,0) {};
\node[fill=white,inner sep=2pt] (dots) at (2,0) {$\cdots$};
\node[dynkin neutral] (b) at (3,0) {};
\draw[dynkin edge] (u)--(c)--(v);
\draw[dynkin edge] (c)--(a)--(dots)--(b);
\end{tikzpicture}}
\newcommand{\ESixDiagram}{%
\begin{tikzpicture}[x=.7cm,y=.7cm,baseline=-.6ex]
\foreach \x/\sty in {0/plus,1/minus,2/plus,3/minus,4/plus}
  \node[dynkin \sty] (v\x) at (\x,0) {};
\node[dynkin minus] (b) at (2,.7) {};
\draw[dynkin edge] (v0)--(v1)--(v2)--(v3)--(v4);
\draw[dynkin edge] (v2)--(b);
\end{tikzpicture}}
\newcommand{\ESevenDiagram}{%
\begin{tikzpicture}[x=.7cm,y=.7cm,baseline=-.6ex]
\foreach \x/\sty in {0/plus,1/minus,2/plus,3/minus,4/plus,5/minus}
  \node[dynkin \sty] (v\x) at (\x,0) {};
\node[dynkin minus] (b) at (2,.7) {};
\draw[dynkin edge] (v0)--(v1)--(v2)--(v3)--(v4)--(v5);
\draw[dynkin edge] (v2)--(b);
\end{tikzpicture}}
\newcommand{\EEightDiagram}{%
\begin{tikzpicture}[x=.7cm,y=.7cm,baseline=-.6ex]
\foreach \x/\sty in {0/plus,1/minus,2/plus,3/minus,4/plus,5/minus,6/plus}
  \node[dynkin \sty] (v\x) at (\x,0) {};
\node[dynkin minus] (b) at (2,.7) {};
\draw[dynkin edge] (v0)--(v1)--(v2)--(v3)--(v4)--(v5)--(v6);
\draw[dynkin edge] (v2)--(b);
\end{tikzpicture}}
\usepackage[hidelinks]{hyperref}
\usepackage{fancyhdr}
\numberwithin{equation}{section}
\setlist[enumerate]{label=\textup{(\roman*)},leftmargin=2em,itemsep=3pt}
\newtheorem{theorem}{Theorem}[section]
\newtheorem{proposition}[theorem]{Proposition}
\newtheorem{lemma}[theorem]{Lemma}
\newtheorem{corollary}[theorem]{Corollary}

\newtheorem{maintheorem}{Theorem}

\theoremstyle{definition}
\newtheorem{definition}[theorem]{Definition}
\newtheorem*{arnoldproblem}{Arnold's problem}

\theoremstyle{remark}\newtheorem{remark}[theorem]{Remark}
\newcommand{\R}{\mathbb R}\newcommand{\C}{\mathbb C}\newcommand{\Z}{\mathbb Z}
\newcommand{\PP}{\mathbb P}\newcommand{\Q}{\mathbb Q}
\newcommand{\HF}{\operatorname{HF}}
\newcommand{\Symp}{\operatorname{Symp}}
\newcommand{\Area}{\operatorname{Area}}
\newcommand{\htop}{h_{\mathrm{top}}}\newcommand{\id}{\operatorname{id}}
\newcommand{\tr}{\operatorname{tr}}
\newcommand{\rk}{\operatorname{rank}}\newcommand{\Int}{\operatorname{int}}
\newcommand{\Fl}{\mathcal B}\newcommand{\Ad}{\operatorname{Ad}}
\newcommand{\Art}{\mathsf A}\newcommand{\Endg}{\mathscr S}
\newcommand{\g}{\mathfrak g}\newcommand{\kalg}{\mathfrak k}\newcommand{\talg}{\mathfrak t}
\newcommand{\halg}{\mathfrak h}

\newcommand{\ZSt}{\mathcal Z}

\newcommand{\Hom}{\operatorname{Hom}}
\newcommand{\Lie}{\operatorname{Lie}}

\newcommand{\pr}{\operatorname{pr}}

\hypersetup{pdftitle={Dynkin plumbing, Dehn twists and entropy production},pdfauthor={},pdfsubject={Dynkin plumbing, symplectic monodromy, and global Floer entropy: main-line version}}
\title{\textbf{Dynkin plumbing, Dehn twists\\and entropy production}}
\author{%
  Wenmin Gong\thanks{%
    School of Mathematical Sciences,
    Beijing Normal University, Beijing, China.
    Email: \texttt{wmgong@bnu.edu.cn}.%
  }
  \and
  Zhijing Wendy Wang\thanks{%
    Department of Mathematics,
    University of Chicago, Chicago, United States of America.
    Email: \texttt{zhijingw@uchicago.edu}.%
  }
  \and
  Jinxin Xue\thanks{%
   New Cornerstone Sciences Laborotary, Department of Mathematics,
    Tsinghua University, Beijing, China.
    Email: \texttt{jxue@tsinghua.edu.cn}.%
  }
}

\date{}
\begin{document}
\maketitle
\begin{abstract}
Exploring Arnold's line of ideas from Hilbert's thirteenth problem through braid groups to symplectic monodromy, we study monodromy of coadjoint orbits, symplectic Dehn twists and the dynamical complexity of their compositions. In his work
on Milnor fibrations, Arnold studied coadjoint orbits and introduced Dehn
twists on cotangent bundles of flag manifolds. We answer two questions
posed in this setting. First, for every connected complex semisimple
group, we construct complete flat symplectic transport on the regular
coadjoint-orbit fibration, matching a prescribed exact real KKS form on
a reference orbit. Second, under an explicit condition at infinity,
the simple-root family twists define a representation of the
ADE Artin group in the equivariant symplectic mapping-class group of
the flag cotangent bundle, which is faithful in type $A$.

We then study entropy production in three geometric settings governed
by ADE diagrams. For sphere plumbings associated with simple hypersurface
singularities, we extend our earlier type-$A$ results to types $D$ and $E$,
obtaining entropy lower bounds in terms of pseudo-Anosov stretch factors
on associated surfaces.  
In the representation-theoretic setting of nilpotent cones and Springer
resolutions, we pass to subregular Slodowy slices, whose ADE surface
fibers carry compactly supported sphere-twist monodromy satisfying the
same bounds. Finally, for ADE configurations of
Lagrangian spheres on a polarized K3 surface, we prove a positive entropy lower bound
and exponential growth of  global Floer ranks of bipartite positive-negative even
twist words. 
\end{abstract}
\setcounter{tocdepth}{1}
\tableofcontents

\section{Introduction}\label{sec:historical-intro}
\label{subsec:main-results}

This paper grows out of an effort to understand an important part of Arnold's stream of mathematical ideas starting from Hilbert's 13 problem, in particular, the paper \cite{ArnoldMilnor}. The basic setting is as follow. Fix $N\ge2$, we consider the set of degree $N$ monic polynomials, normalized so that the roots sum to zero:
\begin{equation}\label{eq:parameter-space}
 \mathcal P_N=\{P_b(z)=z^N+b_2z^{N-2}+\cdots+b_N:
                    b\in\C^{N-1}\}.
\end{equation}
The polynomials with repeated roots form the
\emph{discriminant locus} $\mathcal D_N$ whose complement is the set of regular polynomials
\begin{equation}\label{eq:discriminant-base}
 \mathcal U_N=\mathcal P_N\setminus\mathcal D_N.
\end{equation}
Consider the covering $\{(P_b,z):P_b\in \mathcal U_N,\, P_b(z)=0\}$ over $\mathcal U_N$. Then a loop in $\mathcal U_N$ traces a braid in the lift. Thus we have
\begin{equation}\label{eq:braid-base}
             \pi_1(\mathcal U_N,b_0)=Br_N,
\end{equation}
where $Br_N$ is the braid group on $N$ strands \cite{FadellNeuwirth}.
We fix a base point $P_{b_0}$ with distinct real roots and order them increasingly.

One reason to study the topology of $\mathcal U_N$ comes from Hilbert's thirteenth problem, which asks whether solutions of polynomial equations can be expressed compositions of functions of fewer variables. The superposition theorems of Arnold and Kolmogorov~\cite{Arnold1957,Kolmogorov1957} gives an affirmative answer for compositions of \emph{continuous} functions of two variables. However, a perhaps more natural question is to ask whether solutions can be expressed as composition of \emph{algebraic} functions. Arnold sought obstructions in the topology of algebraic functions and their discriminants, which led to his studies of braid groups and their cohomology \cite{Arnold1969,Arnold1970I,Arnold1970II,ArnoldSuperposition}.
This leads to a geometric realization problem: can we find a (faithful) representation of the braid group as appropriate mapping class group on a manifold? 

In \cite{ArnoldMilnor}, Arnold proposed studying the monodromy of the regular fibers of the adjoint quotient
\begin{equation}\label{eq:full-adjoint-quotient}
 \chi_N:\mathfrak{sl}_N(\C)\longrightarrow\mathcal P_N,
                 \qquad \chi_N(x)=\det(zI-x).
\end{equation}
Over $\mathcal U_N$, each fiber is a single conjugacy orbit under $SL_N(\C)$. For the restriction
\begin{equation}\label{eq:orbit-base}
 \mathcal E_N=\chi_N^{-1}(\mathcal U_N),\qquad
 \pi_N=\chi_N|_{\mathcal E_N}:\mathcal E_N\longrightarrow\mathcal U_N,
 \quad\text{ with fibers }\quad X_b=\pi_N^{-1}(b),
\end{equation}
the regular fibers can be identified, as smooth manifolds, with $T^*\Fl_N$, the cotangent bundle of the flag manifold $\Fl_N$ of $\C^N$. Here the flag manifold $\Fl_N$ is the set of complete flags
\[
 0=V_0\subset V_1\subset\cdots\subset V_N=\C^N,
 \qquad \dim_\C V_j=j
\] and its cotangent bundle $T^*\Fl_N$ can be identified with the Springer resolution of the nilpotent cone \[T^*\Fl=\{(V_\bullet,Z):ZV_j\subset V_{j-1}\}.\]

 Endow the fiber $X_{b_0}$ with the natural symplectic form on cotengent bundle, which is equivalent to the real KKS form $\omega_0$ (see Definition~\ref{def:reference-kks}), Arnold asks whether there exists a symplectic form across the fibers whose monodromy become symplectomorphisms of $X_{b_0}$. He then introduced a symplectic Dehn twist as proposed generators of the monodromy group, whose $T^*\mathbb S^N$ version, the \emph{symplectic Dehn twists}, turns out to be important in sympletic topology after the extensive work of Seidel et al. 

One way of revealing the complexity is from the dynamics, i.e. studying iteration of composite symplectic Dehn twists. In \cite{GWX}, the author studied these compositions of symplectic Dehn twists for $A_n$-type singularities, and used it to get a higher dimensional generalization of the Nielsen-Thurston's classification of mapping class groups on surfaces. We used two indicators for the complexity of the dynamics: topological entropy and growth of ranks of Floer homology group. The first contribution of this paper is to extend the results of \cite{GWX} to $D,E$-type singularities and plumbing of cotangent bundle of spheres for any tree. 

Arnold's Dehn twists on $T^*\Fl_N$ is much less studied. In representation theory, the cotangent bundle of flag manifolds is the total space of the Springer resolution of the nilpotent cone.  Arnold's Dehn twists in this setting are defined  nilpotent elements the Lie group $\mathrm{SL}_N\C$, which in general do not have exponential growing complexity when studying iterations of compositions. However, we show in this paper that dynamical complexity such as postie topological entropy shows up after we taking a Slodowy slice (see Theorem~\ref{main:slice}). Along the line, we also study the basic properties of the fibration and address two open problems of Arnold in \cite{ArnoldMilnor} in Theorems~\ref{main:flat-orbits} and~\ref{main:flag}, which gives a compatible symplectic structure on all regular fibers of $\pi_N$, and give monodromies under different boundary conditions of $T^*\Fl_N$. 

We also give an application of the aforementioned results to the $(-2)$-curves in K3 surfaces, which also admits ADE-configurations. Compactness of the K3 surfaces allows us to study the complexity of the global Floer homology group for composite Dehn twists under iterations. The precise statement is given in Theorem~\ref{main:k3}.  

As we have mentioned above, this paper is motivated by an attempt to understanding an important line of Arnold's research. The results of this paper shed lights on several problems of Arnold asked in  \ref{ArnoldProblems}. We will discuss these problems in Appendix \ref{subsec:arnold-problems}. 
\subsection{Flat transport on regular adjoint orbits}
The following problem was asked in \cite{ArnoldMilnor}.
\begin{arnoldproblem}
Fix a real-spectrum regular orbit $X_{b_0}$ with its prescribed exact real symplectic form $\omega_0$. Can the other regular orbits be equipped with smoothly varying real symplectic forms $\omega_b$, all symplectomorphic to $\omega_0$, so that a based loop $\gamma$ admits symplectic transports
\[
 s_t:(X_{b_0},\omega_0)\longrightarrow(X_{\gamma(t)},\omega_{\gamma(t)}),
 \qquad s_0=\id,
\]
and these transports can be chosen with return maps defining a homomorphism $\pi_1(\mathcal U_N,b_0)\to\Symp(X_{b_0},\omega_0)$?
\end{arnoldproblem}
We first construct a closed real two-form on $\mathcal E_N$, symplectic on each fiber, whose complete flat transports give a homomorphism $\pi_1(\mathcal U_N,b_0)\to\Symp(X_{b_0},\omega_0)$ in Theorem~\ref{main:flat-orbits}.

The construction also works more generally for other complex semisimple Lie groups. Let $G$ be a connected complex semisimple Lie group (examples include $SL_n(\C)$ for $n\ge2$, $Sp_{2n}(\C)$ for $n\ge1$, and $SO_n(\C)$ for $n\ge3$, as well as the complex exceptional groups of types $G_2,F_4,E_6,E_7,E_8$).
Let $K\subset G$ be a maximal compact subgroup and $T\subset K$ a
maximal torus. Write $\kalg=\Lie(K)$ and $\talg=\Lie(T)$. The
complexification $T_\C$ of $T$ has Lie algebra $\halg=\talg_\C$,
a Cartan subalgebra of $\g=\Lie(G)$. A torus is a compact connected abelian Lie group. Let $\Gamma$ be the Dynkin diagram of a chosen system of simple roots. 

\begin{figure}[tbp]
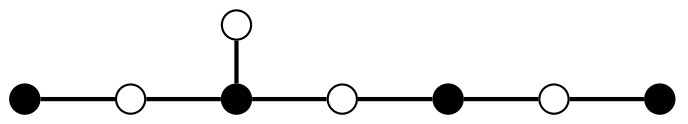

\centering
\renewcommand{\arraystretch}{1.85}
\begin{tabular}{@{}cll@{}}
\toprule
Type & Diagram & Structure / arm lengths\\
\midrule
$A_2$ & \AiiDiagram & two spheres meeting once\\
$A_r$ & \ArDiagram & a chain of $r$ vertices\\
$D_r$ & \DtreeDiagram & $(1,1,r-3)$\\
$E_6$ & \ESixDiagram & $(1,2,2)$\\
$E_7$ & \ESevenDiagram & $(1,2,3)$\\
$E_8$ & \EEightDiagram & $(1,2,4)$\\
\bottomrule
\end{tabular}
\caption{The ADE-Dynkin diagrams. For the flag cotangent bundle, vertices label the simple roots and their family twists. \\In a plumbing or Slodowy fiber, each vertex denotes a sphere and each edge a transverse intersection. 
}
\label{fig:dynkin-main}
\end{figure}
The Weyl group $\mathsf W_\Gamma=N_K(T)/T$ is the finite group generated by the root reflections. Set
\begin{equation}\label{eq:ADE-parameter-space}
 \begin{aligned}
 \halg_{\rm reg}&=\{\lambda\in\halg:
                        \alpha(\lambda)\ne0\text{ for every root }\alpha\},\\
 \mathcal P_\Gamma&=\halg/\mathsf W_\Gamma,\qquad
 \mathcal U_\Gamma=\halg_{\rm reg}/\mathsf W_\Gamma.
 \end{aligned}
\end{equation}
We have the adjoint quotient $\chi_\Gamma:\g\to\mathcal P_\Gamma$ whose regular
part is
\begin{equation}\label{eq:general-orbit-base}
 \mathcal E_\Gamma=\chi_\Gamma^{-1}(\mathcal U_\Gamma),\qquad
 \pi_\Gamma:\mathcal E_\Gamma\longrightarrow\mathcal U_\Gamma, \quad\text{ with fibers }\quad
  X_b=\pi_\Gamma^{-1}(b).
\end{equation}

Each fiber $X_b$ is an adjoint orbit, diffeomorphic to $G/T_\C$ and hence to $T^*\Fl$, where $\Fl=K/T$ is the full flag manifold.
The fundamental group of $\mathcal U_\Gamma$ is the Artin group $\Art_\Gamma$ \cite{BrieskornBraid}, which has generators $\sigma_i$ corresponding to the simple reflections $s_i$ and relations
\[
 \underbrace{\sigma_i\sigma_j\sigma_i\cdots}_{m_{ij}\text{ terms}}
 =\underbrace{\sigma_j\sigma_i\sigma_j\cdots}_{m_{ij}\text{ terms}}
 \quad(i\ne j),
\]
where $m_{ij}\in \{2,3,4,6\}$ is the order of $s_is_j$ in $\mathsf W_\Gamma$.
We have a natural quotient $p:\Art_\Gamma\to\mathsf W_\Gamma;\;\sigma_i\mapsto s_i$, obtained by imposing $\sigma_i^2=1$. For $G=SL_N(\C)$, $p$ is the permutation map $Br_N\to S_N$.

Fix a regular orbit $X_{b_0}$ with real adjoint spectrum. Equip the reference orbit $X_{b_0}=G\cdot\lambda^0$ with the real part
$\omega_0$ of the holomorphic Kirillov--Kostant--Souriau form
(Definition~\ref{def:reference-kks}). Write $\theta_{\rm can}$ for the
tautological one-form on $T^*\Fl$.  Our first result realizes Arnold's
transport on the full family.

\begin{maintheorem}\label{main:flat-orbits}
For every connected complex semisimple $G$ 
there is a smooth real one-form $\Theta$ on $\mathcal E_\Gamma$ such that
\[
 \omega_b=(d\Theta)|_{X_b}\text{ is symplectic},\qquad
 (X_b,\omega_b)\cong(T^*\Fl,d\theta_{\rm can}),\qquad
 \omega_{b_0}=\omega_0.
\]
The distribution $\mathcal H=\ker d\Theta$ is a complete flat connection
for $\pi_\Gamma$, and its parallel transports preserve the fiber forms.
Under a symplectic identification of the reference fiber with
$T^*\Fl$, its holonomy is a homomorphism
\begin{equation}\label{eq:flat-holonomy-intro}
 \begin{aligned}
 \Art_\Gamma&\xrightarrow{\ p\ }\mathsf W_\Gamma
       \xrightarrow{\ w\mapsto\widehat r_w\ }\Symp(T^*\Fl).
 \end{aligned}
\end{equation}
\end{maintheorem}
Arnold identified the real-spectrum orbit in type $A$ symplectically
with a cotangent bundle. Azad, van den Ban, and Biswas extended this
identification to complex semisimple groups and partial flags
\cite{ArnoldMilnor,ABB}. We use the orbitwise identification to match
the form at $b_0$. The closed form on the entire regular family comes
from the product with the Weyl cover. Its null directions give
complete transport, and pure Artin braids, the elements of $\ker p$,
have trivial return maps.

\subsection{Dehn twists and the braid action}\label{subsec:intro-twists}
The construction of Theorem~\ref{main:flat-orbits} answers the question in \cite[p.~101]{ArnoldMilnor}, but its monodromy factors through the Weyl group and is therefore not faithful. The reason is partly that the construction imposes no control at infinity. To retain more information, we choose symplectically isotopic representatives in the same classes with prescribed behavior at infinity.
The boundary condition can be stated directly in terms of the Springer map $\mu_\C:T^*\Fl\to\mathcal N$, where $\mathcal N\subset\g$ is the nilpotent cone. In type $A$, the cotangent bundle is $T^*\Fl=\{(V_\bullet,Z):ZV_j\subset V_{j-1}\}$  identified with flags and $\mu_\C(V_\bullet,Z)=Z$ simply forgets the flag.

Fix a norm on the cotangent bundle, for $p=(x,\zeta)\in T^*\Fl$, write $|p|=|\zeta|$ and $rp=(x,r\zeta)$. A nilpotent element $Z$ is \emph{regular} if $\mu_{\C}^{-1}(Z)$ is a single point.

\begin{definition}[Boundary condition at infinity]\label{def:end}
The group $\Endg_K(T^*\Fl)$ consists of the $K$-equivariant
symplectomorphisms $f$ satisfying:
\begin{enumerate}
\item \emph{The rescaled Springer images agree at infinity:}
\begin{equation}\label{eq:end-rescaled}
 \lim_{R\to\infty}\ \sup_{|p|\ge R}
 \frac{|\mu_\C(f(p))-\mu_\C(p)|}{|p|}=0.
\end{equation}
\item \emph{The map is exactly the identity sufficiently far out in
regular directions:} for every compact set
\[
 C\subset\{p\in T^*\Fl:|p|=1,\ \mu_\C(p)\text{ is regular}\},
\]
there is $R_C>0$ such that $f(rp)=rp$ for all $p\in C$ and $r\ge R_C$.
\end{enumerate}
\end{definition}

Section~\ref{subsec:condition-infinity} specifies the topology on $\Endg_K(T^*\Fl)$. In type $A$ the two conditions above agree with
Jin's definition \cite[Definition~2.1]{Jin}, as we show in
Section~\ref{subsec:braid-representation}.

Following Arnold's construction, we show that there are representatives the monodromy, satisfying the boundary condition of $\Endg_K(T^*\Fl)$, which gives a faithful representation from the Artin group to $\pi_0 \Endg_K(T^*\Fl)$ for type $A$ groups.

\begin{maintheorem}\label{main:flag}
Let $G$ be a connected complex simple Lie group with Dynkin diagram
$\Gamma$ of type $A_r$, $D_r$, $E_6$, $E_7$, or $E_8$. Then there exists a family of Dehn twist $\tau_i$ on $T^*\Fl$, where $i$ ranges over the simple roots of $G$, which define a homomorphism
\[
 \rho_{\mathrm{flag}}:\Art_\Gamma\longrightarrow
           \pi_0\Endg_K(T^*\Fl),\qquad
                 \sigma_i\longmapsto[\tau_i].
\]
Both the maps and the isotopies establishing the Artin relations
satisfy the condition at infinity of Definition~\ref{def:end}.
The homomorphism $\rho_{\mathrm{flag}}$ is faithful when $\Gamma=A_r$.
\end{maintheorem}
At a covector restricting to zero on the sphere fiber, $\tau_i$
acts by the antipodal map along the fiber. For nonzero
restriction, it follows the round geodesic flow with a cutoff in
the norm of that restriction.

Jin's homomorphism for the type-$A$ flag detects all path
components in ranks one and two \cite{Jin}. For adjacent simple
roots we use his rank-two result and extend the resulting path
from $SU(3)/T^2$ to the full flag cotangent bundle.  Here the faithfulness is only established for $\Gamma=A_r$ and  remains an open problem for types $D,E$. 

\subsection{Slodowy slices, sphere twists, and entropy production}
The twists constructed in Theorem~\ref{main:flag} have unbounded support when the rank is at least two. To obtain compactly supported monodromy, we may change the fibers. Brieskorn and Slodowy \cite{Brieskorn,Slodowy} relate the fiber $X_b$ to simple surface singularities by restricting it to a transverse slice called the Slodowy slice. Here we take $G$ to be simple of ADE type.

Let $r=\dim_\C\halg$ be the rank of $G$, and choose a subregular nilpotent element $e\in\g$, meaning that its orbit has complex dimension $\dim_\C\g-r-2$. Complete it to a triple $e,h,f$ satisfying the $\mathfrak{sl}_2\C$ relations $[h,e]=2e$, $[h,f]=-2f$, and $[e,f]=h$. The affine space
\[\mathcal S_e=e+\ker(\operatorname{ad}_f)\subset \mathfrak g\] is the \emph{Slodowy slice}. It is transverse to the orbit of $e$ and has complex dimension $r+2$. We consider the restricted family
\begin{equation}\label{eq:Slodowy-family}
 \chi_{\Gamma,e}=\chi_\Gamma|_{\mathcal S_e}:
                   \mathcal S_e\longrightarrow\mathcal P_\Gamma.
\end{equation}
For $G$ of ADE type, the singular fiber $\chi_{\Gamma,e}^{-1}(0)$ containing $e$ has an isolated surface singularity, and its regular fibers have complex dimension two. 

Choose a sufficiently small closed ball $\mathcal K_e\subset\mathcal S_e$ around $e$, called a \emph{Milnor neighborhood}. For a sufficiently small
neighborhood $0\in \mathcal V_\Gamma\subset \chi_{\Gamma,e}(\mathcal K_e)$ in the base, let
$\mathcal U_\Gamma^{\rm loc}=\mathcal V_\Gamma\cap\mathcal U_\Gamma$.
Restricting the regular fibers to these neighborhoods gives
\begin{equation}\label{eq:compact-slice-fibration}
 \pi_\Gamma^{\rm slice}:\mathcal F_\Gamma
       =\mathcal K_e\cap\chi_{\Gamma,e}^{-1}(\mathcal U_\Gamma^{\rm loc})
       \longrightarrow\mathcal U_\Gamma^{\rm loc},\qquad
 F_{\Gamma,b}=\mathcal K_e\cap\chi_{\Gamma,e}^{-1}(b).
\end{equation}
The neighborhood can be chosen so that
$\mathcal U_\Gamma^{\rm loc}\hookrightarrow\mathcal U_\Gamma$
induces an isomorphism on fundamental groups; see Section~\ref{sec:slodowy}. Each fiber $F_{\Gamma,b}$ is a compact four-manifold with boundary, called a \emph{Milnor fiber}, with its standard exact symplectic structure. Write $F_\Gamma=F_{\Gamma,b_*}$ for a fixed point $b_*\in \mathcal U_\Gamma^{\rm loc}$.
Using a boundary collar trivialization over $\mathcal V_\Gamma$, transports normalized to be the identity near the boundary give compactly supported monodromies on $\Int F_\Gamma$.

The Milnor fiber is diffeomorphic to a plumbing of disk cotangent bundles $D^*S^2$ according to the graph $\Gamma$ (see Section~\ref{sec:tree}). There is one Lagrangian sphere $L_v$ for each vertex of $\Gamma$; two spheres intersect transversely once if their vertices are joined by an edge and are otherwise disjoint. The Lagrangian spheres $L_v$ are vanishing spheres that collapse along suitable paths when moving into fibers of the discriminant. Symplectic Picard--Lefschetz theory identifies monodromy around simple discriminant points with Dehn twists about the corresponding spheres \cite{SeidelLES,SeidelBook}. This gives a homomorphism
\[
 \rho_{\mathrm{slice}}:\Art_\Gamma\longrightarrow
         \pi_0\Symp_c(\Int F_\Gamma),\qquad
              \sigma_v\longmapsto[\tau_{L_v}]
\]
where the standard Artin generators are represented by the symplectic Dehn twists.

This construction restricts the deformation family and then takes monodromy; the representatives in Theorem~\ref{main:flag} need not preserve the slices.  The  representation $\rho_{\mathrm{slice}}$ is faithful in every ADE type by Brav--Thomas's faithfulness theorem \cite{BravThomas} and by Khovanov--Seidel \cite{KhovanovSeidel} for type $A$. 

Here, we discuss a further question: how much dynamical complexity is forced by the symplectic mapping classes? We measure complexity by \emph{topological entropy}, denoted $\htop$, which measures the exponential growth rate of distinguishable orbit segments.

Our construction works for any finite tree $\Gamma$ (not necessarily a Dynkin diagram). 
For any finite tree $\Gamma$, let $W_\Gamma$ be the plumbing of $D^*S^2$ along $\Gamma$, with vertices $v$ representing Lagrangian spheres $L_v$ and edges representing intersection, and let $\tau_{L_v}$ be the symplectic Dehn twist around $L_v$. We also take $\Sigma_\Gamma$ to be the surface obtained by plumbing of annuli $D^*S^1$ along $\Gamma$, with a circle $\gamma_v$ for each vertex $v$, and let $T_{\gamma_v}$ be the Dehn twist around $\gamma_v$ in $\Sigma_\Gamma $.The Artin group $\Art_\Gamma$ has generators $\sigma_v$, with the three-term braid relation on edges and commutation for nonadjacent vertices. 

For a word $\beta=\sigma_{v_\ell}^{e_\ell}\cdots\sigma_{v_1}^{e_1}
\in\Art_\Gamma$, $e_j\in\Z$, take
\begin{equation}\label{eq:arbitrary-tree-word}
\Phi_\beta=\tau_{L_{v_\ell}}^{e_\ell}\cdots\tau_{L_{v_1}}^{e_1},\qquad
 \phi_\beta=T_{\gamma_{v_\ell}}^{e_\ell}\cdots T_{\gamma_{v_1}}^{e_1},\qquad
 M_\beta=(\phi_\beta)_*|_{H_1(\Sigma_\Gamma;\R)}.
\end{equation}
We also denote by $\Lambda_\beta$ the maximum of $1$ and the largest stretch factor of $\phi_\beta$ among its pseudo-Anosov components. See Definition~\ref{def:pA-growth} for the precise normalization.

Write $\Symp_c(\Int W_\Gamma)$ for smooth symplectomorphisms equal to the identity near the boundary. Then we have the following bound.

\begin{maintheorem}\label{main:tree}
For any finite tree $\Gamma$ and any word $\beta\in\Art_\Gamma$, we have
\begin{equation}\label{eq:tree-class-entropy}
 \inf_{\substack{f\in\Symp_c(\Int W_\Gamma)\\{}[f]=[\Phi_\beta]}}
       \htop(f)
       \ge\log\Lambda_\beta\ge\log\rho(M_\beta).
\end{equation}
Here $\rho(M_\beta)$ is the spectral radius.
\end{maintheorem}

In \cite{GWX} we proved positive entropy for selected products on
chains of spheres and used Floer ranks to study their iterates.
Keating's Lefschetz stabilization compares the Floer ranks of
arbitrary signed twist words in a plumbing with those of the
associated surface word \cite{KeatingStab}. We use the
Nielsen--Thurston decomposition to measure the surface growth, then
pass from Floer ranks to area growth and entropy
\cite{FLP,CGG,Yomdin}. 

Applying the plumbing estimate to the slice monodromy gives the
following result.

\begin{maintheorem}\label{main:slice}
For every ADE diagram $\Gamma$, the boundary-normalized monodromy of
\eqref{eq:compact-slice-fibration} gives a faithful homomorphism
\[
 \rho_{\mathrm{slice}}:\Art_\Gamma\longrightarrow
       \pi_0\Symp_c(\Int F_\Gamma),\qquad
       \sigma_v\longmapsto[\tau_{L_v}].
\]
For every $\beta\in\Art_\Gamma$ and every representative
$f\in\Symp_c(\Int F_\Gamma)$ of $\rho_{\mathrm{slice}}(\beta)$,
\[
              \htop(f)\ge\log\Lambda_\beta\ge\log\rho(M_\beta).
\]
\end{maintheorem}

For type $A$, we use Khovanov--Seidel's monodromy theorem, whose proof
relates Floer ranks to intersections of arcs in a marked disk
\cite{KhovanovSeidel}. Appendix~\ref{app:slice-faithful-de} proves
faithfulness for types $D$ and $E$ using Brav--Thomas's theorem on
configurations of spherical twists \cite{BravThomas}. The
dynamical estimate applies in all three types: if
$\Lambda_\beta>1$, every representative of the monodromy has positive
entropy.

\subsection{Floer homology growth and positive entropy on K3 surfaces}
An ADE configuration of Lagrangian spheres in a closed symplectic
manifold has a plumbing neighborhood. Extending its twists by the
identity gives global symplectomorphisms. To bound entropy throughout
the ambient symplectic class, we compute Floer groups in the entire
closed manifold. Our final result treats bipartite positive-negative
even twist words on polarized K3 surfaces.

A \emph{K3 surface} is a compact simply connected complex surface with
trivial canonical bundle.
\begin{definition}[Polarized K3]\label{def:polarized-k3}
A \emph{polarization} of a projective K3 surface $(X,I_0)$ is an ample
line bundle $\mathcal H$. We assume that $\omega$ is a K\"ahler form
for $I_0$ and $[\omega]=c\,c_1(\mathcal H)$ for some $c>0$.
We call $(X,I_0,\omega)$ a symplectic polarized K3.
\end{definition}

Let $\Gamma$ be a connected ADE diagram with at least two vertices.
An ADE configuration consists of embedded Lagrangian two-spheres,
one for each vertex, with one transverse intersection for each edge
and no other intersections. Label its two colors
$A_1,\ldots,A_m$ and $B_1,\ldots,B_\ell$, and set
\[
 C_{ij}=\begin{cases}1,&A_i\text{ meets }B_j,\\0,&\text{otherwise},\end{cases}
 \qquad \sigma=\|C\|_2>0.
\]
Spheres of one color are disjoint, so their twists commute.
For integers $p,q\ge1$, put
\begin{equation}\label{eq:K3-word}
 \Phi_{p,q}=
 \left(\prod_{i=1}^m\tau_{A_i}^{2p}\right)
 \left(\prod_{j=1}^{\ell}\tau_{B_j}^{-2q}\right),
\end{equation}
\begin{equation}\label{eq:K3-lambda}
 \lambda_\Gamma(p,q)=
 1+2pq\sigma^2+2\sqrt{pq\sigma^2(1+pq\sigma^2)}.
\end{equation}
Composition is from right to left.

\begin{maintheorem}\label{main:k3}
Suppose a symplectic polarized K3 $(X,I_0,\omega)$ contains the above
ADE configuration. There is a fixed twist descendant $K_*$ of a core
sphere such that, for every core $L_v$,
\begin{equation}\label{eq:K3-growth}
 \lim_{n\to\infty}\frac1n
 \log\dim_\Lambda\HF_X(L_v,\Phi_{p,q}^{\,n}K_*)
       =\log\lambda_\Gamma(p,q).
\end{equation}
Here $\Lambda$ is the universal Novikov field with complex coefficients,
defined at the start of Section~\ref{sec:k3-main}. All Floer groups
are computed in $X$.
For every nonzero integer $k$,
\begin{equation}\label{eq:entropy}
 \inf_{\substack{f\in\Symp(X,\omega)\\{}[f]=[\Phi_{p,q}^{\,k}]}}
 \htop(f)\ge |k|\log\lambda_\Gamma(p,q)>0.
\end{equation}
Each $\Phi_{p,q}^{\,k}$ is smoothly isotopic to the identity and acts
trivially on integral cohomology.
\end{maintheorem}
Seidel showed that the square of a sphere twist on a K3 can be
smoothly trivial and symplectically nontrivial. Sheridan and Smith
later showed that the symplectic Torelli group, the subgroup of
mapping classes acting trivially on integral cohomology, can be
infinitely generated \cite{SeidelK3Twist,SheridanSmithK3}.
Theorem~\ref{main:k3} gives an entropy estimate for explicit words
in the kernel of the map from symplectic to smooth mapping classes.
The main step is to control the product maps in Seidel's closed
twist triangle \cite{SeidelQuartic}. Their vanishing gives a linear
recurrence for the ranks, even at a vertex with several neighbors.
Its Perron eigenvalue is $\lambda_\Gamma(p,q)$. For $D_4$ with the
central sphere taken as the first color,
\[
 C=(1\ 1\ 1),\qquad
 p=q=1\quad\Longrightarrow\quad \lambda_{D_4}(1,1)=7+4\sqrt3 .
\]
For an $A_2$-configuration, one has $\sigma=1$ and we get the following results.

\begin{corollary}\label{cor:k3-A2}
Let $L_1,L_2$ be Lagrangian spheres meeting transversely once on a
symplectic polarized K3. For $p,q\ge1$, define
\begin{equation}\label{eq:Psi}
 \Psi_{p,q}=\tau_{L_1}^{2p}\tau_{L_2}^{-2q},
\end{equation}
\begin{equation}\label{eq:lambda}
 \lambda_{p,q}=1+2pq+2\sqrt{pq(pq+1)},\qquad
 D_n(p,q)=\frac{\lambda_{p,q}^{\,n}+\lambda_{p,q}^{\,1-n}}
                   {\lambda_{p,q}+1}.
\end{equation}
For every $n\ge0$,
\begin{equation}\label{eq:growth}
 \dim_\Lambda\HF_X(L_1,\Psi_{p,q}^{\,n}L_2)=D_n(p,q).
\end{equation}
\end{corollary}
For $p=q=1$, these ranks begin as $1,1,5,29,169,985$, with exponential
rate $\log(3+2\sqrt2)$.

An $A_2$ configuration occurs, for example, in the quartic family
\[
 X_t=\{x_0^2x_3^2+x_1^2x_3^2+x_2^3x_3
              +x_0^4+x_1^4+x_2^4=t x_3^4\}\subset\PP^3.
\]
In the chart $x_3=1$, the defining function is
$F=x^2+x^4+y^2+y^4+z^3+z^4$. Its critical values are precisely the
sums of one entry from each of $\{0,-1/4\}$, $\{0,-1/4\}$ and
$\{0,-27/256\}$. Only the origin has critical value $0$.
At $x_3=0$, the three derivatives $4x_j^3$, $j=0,1,2$, cannot
vanish together. Near the origin the change of coordinates
\[
 u=x(1+x^2)^{1/2},\quad v=y(1+y^2)^{1/2},\quad
 w=z(1+z)^{1/3},\qquad F=u^2+v^2+w^3
\]
exhibits the $A_2$ singularity. Thus for sufficiently small $t\ne0$
the smooth quartic $X_t$ is a K3 with two vanishing spheres meeting
once, Lagrangian for the restricted Fubini--Study form
\cite{SeidelBook}.

The paper is organized as follows. Sections~\ref{sec:flat-orbits} and~\ref{sec:flag} treat the full
orbits and flags. Section~\ref{sec:tree} proves the entropy estimate
for tree plumbings, and Section~\ref{sec:slodowy} applies it to
transverse monodromy. Section~\ref{sec:k3-main} proves the global
Floer growth and the K3 entropy bound.
Appendix~\ref{app:seidel-comparison} gives the closed twist sequence
used there. Appendix~\ref{subsec:arnold-problems} returns to Arnold's
questions, and Appendix~\ref{app:slice-faithful-de} proves faithfulness
of the slice monodromy in types $D$ and $E$.

\section*{Acknowledgment}
J. X. is supported by the National Natural Science Foundation of China
(Grant Nos.~12271285 and 11790273), the New Cornerstone Investigator
program, and the Xiaomi Foundation. J. X. is a member of the IAS
this semester, supported by the Shiing-Shen Chern Member Fund.

\subsection*{AI declaration}
We began this project three years ago, after our work \cite{GWX}. We obtained
most of the results in Theorems A--D before using AI tools. Theorem E
grew out of repeated discussions with AI. We also used AI tools
for exposition and proof checking. %The authors are responsible for the mathematical content.

\section{Flat transport on regular adjoint orbits}\label{sec:flat-orbits}

Let $G$ be a connected complex semisimple Lie group, $K$ a maximal
compact subgroup, and $T$ a maximal torus of $K$. For
$w\in\mathsf W_\Gamma=N_K(T)/T$, choose a representative $\dot w\in N_K(T)$.
Let $r_w(kT)=k\dot w^{-1}T$ on $\Fl=K/T$ and denote its cotangent lift
by $\widehat r_w$. Fix an invariant complex bilinear form $\kappa$ on
$\g$ equal to a positive multiple of the Killing form on each simple
factor. The inner product $\langle\,,\,\rangle=-\kappa|_{\kalg}$
identifies $T^*\Fl$ with $K\times_T\talg^\perp$. In these coordinates,
\[
 \widehat r_w[k,\xi]=[k\dot w^{-1},\Ad_{\dot w}\xi].
\]
The relation $[kt,\Ad_{t^{-1}}\xi]=[k,\xi]$ for $t\in T$ makes the
formula independent of the representative. This gives a holomomy of the bundle $\mathcal E_\Gamma\to \mathcal U_\Gamma$.

\begin{lemma}\label{lem:orbit-product}
There is a bundle diffeomorphism over $\mathcal U_\Gamma$
\begin{equation}\label{eq:orbit-cotangent-quotient}
 (T^*\Fl\times\halg_{\rm reg})/\mathsf W_\Gamma
    \xrightarrow{\ \sim\ }\mathcal E_\Gamma,\qquad
 w(x,\lambda)=(\widehat r_w x,w\lambda).
\end{equation}
The product connection descends to a complete flat connection on the
quotient. In the coordinate $x\mapsto[x,\lambda^0]$ its holonomy along
a based loop $\beta$ is $\widehat r_{p(\beta)}$.
\end{lemma}
\begin{proof}
Mostow decomposition \cite[Proposition~2.20]{Mayrand} identifies
$[k,\xi]$ with $gT_\C=k\exp(i\xi)T_\C$ and carries $\widehat r_w$ to
right multiplication by $\dot w^{-1}$. The orbit map
$(gT_\C,\lambda)\mapsto\Ad_g\lambda$ therefore descends to the
quotient. It is bijective there by regular semisimple conjugacy, and
its differential at $(eT_\C,\lambda)$ is the isomorphism
$(\bar U,v)\mapsto[U,\lambda]+v$ from $(\g/\halg)\oplus\halg$ to $\g$.
Upstairs, parallel transport keeps $x$ fixed along every lifted path.
With right-to-left composition, $\beta$ lifts from $\lambda^0$ to
$p(\beta)^{-1}\lambda^0$, and
$[x,p(\beta)^{-1}\lambda^0]=[\widehat r_{p(\beta)}x,\lambda^0]$.
\end{proof}

Let $q$ be the quotient projection in~\eqref{eq:orbit-cotangent-quotient}.
For a choice of bundle identification $J$ from this quotient onto
$\mathcal E_\Gamma$, define the one-form $\Theta$ on $\mathcal E_\Gamma$ by
\[
 (J\circ q)^*\Theta=\pr_1^*\theta_{\rm can}.
\]
Indeed, $\theta_{\rm can}=\langle\xi,k^{-1}dk\rangle$, and
$\widehat r_w$ acts by $\Ad_{\dot w}$ on both entries of the pairing.
Thus $\pr_1^*\theta_{\rm can}$ is invariant and descends.
Upstairs its differential restricts to $d\theta_{\rm can}$ on each
fiber and has kernel tangent to $\halg_{\rm reg}$.  We choose $J$ below so that $d\Theta|_{X_{b_0}}=\omega_0$.

\begin{definition}[Kirillov--Kostant--Souriau form]\label{def:reference-kks}
On an adjoint orbit, the holomorphic Kirillov--Kostant--Souriau
\textup{(KKS)} form $\Omega$ is given by
\begin{equation}\label{eq:orbit-reference-form}
 \Omega_x([U,x],[V,x])=\kappa(x,[U,V]),\qquad U,V\in\g.
\end{equation}
On the reference orbit $X_{b_0}=G\cdot\lambda^0$, set
$\omega_0=\operatorname{Re}\Omega$.
\end{definition}

\begin{proof}[Proof of Theorem~\ref{main:flat-orbits}]
It remains to choose $J$ so that $d\Theta|_{X_{b_0}}=\omega_0$.
Write $m(x)$ for the image of $[x,\lambda^0]$ under the Mostow
identification in Lemma~\ref{lem:orbit-product}.
Since $\operatorname{ad}_{\lambda^0}$ has real eigenvalues,
Azad, van den Ban and Biswas \cite[Theorem~6.1 and Corollary~6.8]{ABB}
give, after precomposing their identification with fiber negation to convert $-d\theta_{\rm can}$ to $d\theta_{\rm can}$, a $K$-equivariant symplectomorphism
$\Phi:(T^*\Fl,d\theta_{\rm can})\to(X_{b_0},\omega_0)$ with
$\Phi[k,0]=m[k,0]=\Ad_k\lambda^0$.

The map $h=m^{-1}\Phi$ fixes the zero section pointwise. For fiber
dilation $\delta_t(q,p)=(q,tp)$, the diffeomorphisms
$\delta_t^{-1}h\delta_t$ extend smoothly at $t=0$ to the normal
derivative of $h$. This $K$-equivariant bundle automorphism acts at
$eT$ by a nonzero complex scalar on each root plane, since the positive
roots give pairwise inequivalent rotation representations of $T$.
Joining each scalar to $1$ thus gives an isotopy $h_s$ from $\id$ to $h$.
In the local product chart through $\lambda^0$, precompose the Mostow
identification on each fiber with $h_{\chi(b)}$, where $\chi$ is a
smooth cutoff supported in the chart and equal to $1$ near $b_0$.
The resulting $J$ satisfies $J[x,\lambda^0]=m h(x)=\Phi(x)$. Hence
\[
 \Phi^*(d\Theta|_{X_{b_0}})=d\theta_{\rm can}=\Phi^*\omega_0.
\]
Lemma~\ref{lem:orbit-product} gives the asserted holonomy in the
symplectic coordinate $\Phi$.
\end{proof}

\section{ADE flag cotangent bundles and Artin relations}\label{sec:flag}

For each simple root $\alpha_i$, choose a homomorphism $\iota_i:SU(2)\to K$
taking the diagonal circle into $T$ and the upper root space of
$\mathfrak{sl}_2(\C)$ onto $\g_{\alpha_i}$ under the complexified
differential.
In the standard type-$A$ model, the subgroup for $\alpha_i=e_i-e_{i+1}$
acts on the plane spanned by the $i$th and $(i+1)$st coordinate vectors
and fixes its complement. For a set $I$ of simple-root indices, let $H_I$
be the subgroup generated by $T$ and the images of $\iota_i$ for $i\in I$.
Write $H_i=H_{\{i\}}$ and $H_{ij}=H_{\{i,j\}}$.

The projection $\Fl=K/T\to K/H_I$ has fiber $H_I/T$. For one or two
vertices, the fiber is given below, with $i\sim j$ denoting adjacency
in the Dynkin diagram.
\begin{equation}\label{eq:small-flag-fibers}
\begin{array}{c|ccc}
 I&\{i\}&\{i,j\},\ i\not\sim j&\{i,j\},\ i\sim j\\ \hline
 H_I/T&S^2&S^2\times S^2&SU(3)/T^2 .
\end{array}
\end{equation}

\subsection{From a sphere twist to a flag twist}\label{subsec:explicit-flag-twist}

We equip each cotangent bundle with the differential of its tautological
one-form.
\begin{lemma}\label{lem:induction}
For $H=H_I$, every $H$-equivariant symplectomorphism $f$ of $T^*(H/T)$
induces a $K$-equivariant symplectomorphism $\operatorname{Ind}_H^K f$
of $T^*\Fl$. This construction respects composition and continuous
paths in the compact-open $C^\infty$ topology. If $I\subset I'$,
then $\operatorname{Ind}_{H_I}^K f
=\operatorname{Ind}_{H_{I'}}^K(\operatorname{Ind}_{H_I}^{H_{I'}} f)$.
\end{lemma}
\begin{proof}
The free $H$-action $a\cdot(k,x)=(ka^{-1},ax)$ on $K\times(H/T)$
has quotient $K/T$ under $(k,hT)\mapsto khT$. The moment map $J$
of its cotangent lift evaluates a covector on the infinitesimal
$H$-action. Thus $J^{-1}(0)$ consists of the covectors vanishing on
the $H$-orbits. Each is the pullback of a unique covector on the
quotient, giving $J^{-1}(0)/H\cong T^*(K/T)$. Under this identification
the tautological form on the product descends to that on $T^*(K/T)$.

For the cotangent action of $H$ on $T^*(H/T)$, the moment map in
compact coordinates is $\mu_H[h,q]=\Ad_h q$. Since $f$ is symplectic
and $H$-equivariant, $\mu_H\circ f-\mu_H$ is a constant invariant
covector. It vanishes on the derived algebra by invariance and on the
center because the center acts trivially on $H/T$. Thus $\id\times f$
preserves $J$ and descends to the required symplectomorphism.
Descent commutes with composition and is continuous on paths in local
quotient charts.

To check induction in stages, write $\xi=q+v$ with
$q\in\Lie(H)\cap\talg^\perp$ and $v\perp\Lie(H)$.
If $f[1,q]=[g,q']$, preservation of $\mu_H$ gives $\Ad_g q'=q$.
The new base point is $kgT$, and the compact $K$-moment remains
$\Ad_k\xi$, since $\id\times f$ fixes the $T^*K$ factor. Hence
\begin{equation}\label{eq:induction-coordinates}
 (\operatorname{Ind}_H^K f)[k,q+v]=[kg,q'+\Ad_{g^{-1}}v].
\end{equation}
For $H_I\subset H_{I'}$, split $v$ into its components in $\Lie(H_{I'})$
and $\Lie(H_{I'})^\perp$. Inducing in two stages applies
$\Ad_{g^{-1}}$ to each component, giving the same formula.
\end{proof}

On $T^*(H_i/T)\cong T^*S^2$, we take the sphere Dehn twist about the zero
section \cite{ArnoldMilnor}, \cite[Lemma~1.8 and (1.7)]{SeidelLES}.
Let $g^t$ be the
unit-speed geodesic flow on the cotangent bundle of the round unit
sphere away from its zero section. Choose $0<r_0<R$ and a smooth
nonincreasing function $\vartheta:[0,\infty)\to[0,\pi/2]$ equal to
$\pi/2$ on $[0,r_0]$ and to zero on $[R,\infty)$. With $q\in S^2$ and
$p\in T_q^*S^2$, the twist is
\begin{equation}\label{eq:rankone-twist}
 t_\vartheta(q,p)=g^{2\vartheta(|p|)}(q,p)\quad(p\ne0),\qquad
 t_\vartheta(q,0)=(-q,0).
\end{equation}
Near the zero section this is the cotangent lift of the antipodal
map. Away from the zero section it is the time-one map of
$2\int_0^{|p|}\vartheta(s)\,ds$, with the convention
$\iota_{X_H}d\theta_{\rm can}=-dH$.

The round metric is $H_i$-invariant, so $t_\vartheta$ is
$H_i$-equivariant. Define $\tau_i=\operatorname{Ind}_{H_i}^K t_\vartheta$.
For the flag projection $\pr_i:\Fl\to K/H_i$, let
$\mathcal C_i\subset T^*\Fl$ consist of the covectors whose restriction
to the sphere fiber is zero. Explicitly, if $F=\pr_i^{-1}(b)$ and
$x\in F$, then $(x,\xi)\in\mathcal C_i$ precisely when
$\xi|_{T_xF}=0$, or equivalently
$\xi=(d\pr_i)_x^*\eta$ for a unique $\eta\in T_b^*(K/H_i)$.

\begin{lemma}\label{lem:explicit-twist}
On $\mathcal C_i$, the twist $\tau_i$ sends $x$ to its antipode $x'$
in $F$. For each $\eta\in T_b^*(K/H_i)$ its action is given by
\[
 \tau_i\bigl(x,(d\pr_i)_x^*\eta\bigr)
       =\bigl(x',(d\pr_i)_{x'}^*\eta\bigr).
\]
It fixes $(x,\xi)$ whenever $\lVert\xi|_{T_xF}\rVert\ge R$.
Its support is unbounded if $\operatorname{rank}G\ge2$.
\end{lemma}
\begin{proof}
We normalize the invariant inner product on $\kalg$ so
that its restriction to each simple-root $\mathfrak{su}(2)$ is
$-\operatorname{Re}\tr(UV)$. The roots of $\Gamma$ have the same
length, so a single normalization suffices. For $[k,\xi]\in T^*\Fl$, write the component of $\xi$ in
$\Lie(H_i)\cap\talg^\perp$ as
\[
 q_i=d\iota_i Q(z_i),\qquad
 Q(z)=\begin{pmatrix}0&z\\-\bar z&0\end{pmatrix},\qquad r_i=|z_i|.
\]
At the base point of the round sphere, $Q(z_i)$ generates a tangent
vector of length $2r_i$. The covector $q_i$ is dual to half this
vector, since polarizing $-\operatorname{Re}\tr(Q(z_i)^2)=2r_i^2$
gives its pairing with every tangent vector. Thus the restriction
of $\xi$ to the fiber has norm $r_i$, and, for $r_i>0$,
$\exp(tQ(z_i)/r_i)$ moves its base point along the corresponding
geodesic through angle $2t$. By~\eqref{eq:rankone-twist}
and~\eqref{eq:induction-coordinates}, with $n_i=\iota_i(Q(1))$,
\begin{equation}\label{eq:familytwist}
 \tau_i[k,\xi]=[kg,\Ad_{g^{-1}}\xi],\qquad
 g=\begin{cases}
 \iota_i\bigl(\exp(\vartheta(r_i)Q(z_i)/r_i)\bigr),&r_i>0,\\
 n_i,&r_i=0 .
 \end{cases}
\end{equation}
For $0<r_i\le r_0$ the exponential is $Q(z_i/r_i)$, which differs
from $Q(1)$ by right multiplication by a diagonal matrix in $SU(2)$.
Thus for $r_i\le r_0$ the twist equals the Weyl cotangent lift
\begin{equation}\label{eq:zero-extension}
 [k,\xi]\longmapsto[kn_i,\Ad_{n_i^{-1}}\xi].
\end{equation}
The underlying map $r_{s_i}$ is antipodal on $F$ and satisfies
$\pr_i\circ r_{s_i}=\pr_i$. Its cotangent lift therefore carries
$(d\pr_i)_x^*\eta$ to $(d\pr_i)_{x'}^*\eta$, proving the formula
on $\mathcal C_i$. For $r_i\ge R$ the profile vanishes.
In rank at least two, choose $0\ne v\perp\Lie(H_i)$. Every point
$[1,tv]$ on the ray $t\ge0$ moves to the base point $n_iT\ne T$,
so the support is unbounded.
\end{proof}

\subsection{The condition at infinity}\label{subsec:condition-infinity}

Recall the two conditions at infinity from Definition~\ref{def:end}.
For $H=H_I$, use the Springer map $\mu_{\C,H}$ of $T^*(H/T)$ for the
positive roots inherited from $G$, and the norm induced from $\kalg$.
When $H=K$, this is the map $\mu_\C$ on $T^*\Fl$.
An $H$-equivariant symplectomorphism $f$ belongs to $\Endg_H(T^*(H/T))$
if it satisfies both conditions below.
\begin{enumerate}
\item The ratio $|\mu_{\C,H}(f(p))-\mu_{\C,H}(p)|/|p|$ tends uniformly
to zero as $|p|\to\infty$.
\item For each compact set $C$ of unit covectors with regular Springer
image, there is an $R_C$ such that $f(rp)=rp$ for $p\in C$ and
$r\ge R_C$. Here regularity is taken in the semisimple part of
$\Lie(H)_\C$.
\end{enumerate}

\begin{definition}[Paths in $\Endg_H$]
A family $(f_s)_{s\in[0,1]}$ in $\Endg_H(T^*(H/T))$ is a continuous
path if it is continuous in the compact-open $C^\infty$ topology and
\begin{equation}\label{eq:end-family}
 \lim_{R\to\infty}\sup_{s\in[0,1]} D_R(f_s)=0,\qquad
 D_R(f)=\sup_{|p|\ge R}
 \frac{|\mu_{\C,H}(f(p))-\mu_{\C,H}(p)|}{|p|}.
\end{equation}
The cutoff $R_C$ in the second condition may depend on $s$.
For a sequence $f_n\to f$ the analogous conditions are compact-open
$C^\infty$ convergence and~\eqref{eq:end-family} with $\sup_n$ in
place of $\sup_s$.
\end{definition}

\begin{lemma}\label{lem:induction-end}
For $H=H_I$, induction sends $\Endg_H(T^*(H/T))$ into
$\Endg_K(T^*\Fl)$ and takes continuous paths to continuous paths
in the sense of~\eqref{eq:end-family}.
\end{lemma}
\begin{proof}
If $I=\varnothing$, then $H=T$ and induction gives the identity on
$T^*\Fl$. Assume $I\ne\varnothing$ and write $\xi=q+v$ as
in~\eqref{eq:induction-coordinates}, with $q\in\Lie(H)\cap\talg^\perp$
and $v\perp\Lie(H)$. In compact coordinates
$\mu_\C[k,\xi]=\Ad_k\nu(\xi)$, where $\nu(\xi)=\mu_\C[1,\xi]$ is the
normalized projection to the positive root spaces. Its restriction gives
$\mu_{\C,H}[h,q]=\Ad_h\nu(q)$.

We first check condition (i). Adding a root in the
subsystem generated by $I$ leaves the simple-root coefficients outside
$I$ unchanged. Thus the positive and negative root spaces outside the
subsystem form separate $H$-invariant sums, and $\nu$ is $H$-equivariant
on $\Lie(H)^\perp$. If $f[1,q]=[g,q']$, the $v$-terms in
\eqref{eq:induction-coordinates} cancel since
$\Ad_g\nu(\Ad_{g^{-1}}v)=\nu(v)$. Hence
\begin{equation}\label{eq:Levi-error}
 \mu_\C\bigl((\operatorname{Ind}_H^K f)[k,\xi]\bigr)-\mu_\C[k,\xi]
 =\Ad_k\bigl(\mu_{\C,H}(f[1,q])-\mu_{\C,H}[1,q]\bigr).
\end{equation}
For a path $f_s$ and any $\varepsilon>0$, condition~\eqref{eq:end-family}
bounds the norm on the right by $\varepsilon|q|$ for all sufficiently
large $|q|$, uniformly in $s$. On the remaining compact set it is
bounded by a constant $C_\varepsilon$ independent of $s$. Since
$|q|\le|\xi|$, division by $|\xi|$ gives~\eqref{eq:end-family} for the
induced path. A constant path gives the assertion for a single map.

We now check condition (ii). Let $C$ be a compact set of unit covectors
$[k,\xi]$ with regular Springer image for $K$. By Kostant's criterion
\cite[Theorem~5.3]{KostantPrincipal}, $\nu(\xi)$ has nonzero component
in every simple-root space. The projection to $\Lie(H)$ retains the
components indexed by $I$, so $\nu(q)$ is regular for $H$ and $q\ne0$.
As $[k,\xi]$ ranges over $C$, all its representatives form a compact
set. Hence $|q|\ge c>0$, and the directions $[1,q/|q|]$ form a compact
set of regular unit covectors for $H$. If $A$ is their cutoff for $f$,
then $r\ge A/c$ implies $f[1,rq]=[1,rq]$. Formula
\eqref{eq:induction-coordinates} now gives
$(\operatorname{Ind}_H^K f)[k,r\xi]=[k,r\xi]$.
This applies separately to each $f_s$, and Lemma~\ref{lem:induction}
gives compact-open $C^\infty$ continuity of the induced path.
\end{proof}

\subsection{Artin relations in types \texorpdfstring{$A$, $D$, and $E$}{A, D, and E}}
\label{subsec:braid-representation}

We use Jin's result~\cite{Jin} and state it in the notation of Definition~\ref{def:end}.
Write $Br_N$ for the braid group on $N$ strands and $\sigma_i$ for
the standard half-twist exchanging strands $i$ and $i+1$, with the
orientation used by Jin.

\begin{theorem}[{Jin \cite[Theorem~1.4, \S\,3.2, and proof of Corollary~3.6]{Jin}}]\label{thm:jin-braid-map}
For $K=SU(N)$ there is a surjective homomorphism
\[
 \beta_K:\Endg_K(T^*\Fl)\longrightarrow Br_N,\qquad
 \beta_K(\tau_i)=\sigma_i\quad(1\le i<N).
\]
For $N=2,3$, $\beta_K$ is a homotopy
equivalence to the discrete group $Br_N$.
\end{theorem}

Jin formulates the boundary condition using the signed graph
$L_f=\{(-p,f(p)):p\in T^*\Fl\}$ of a $K$-equivariant
symplectomorphism $f$ \cite[Definition~2.1]{Jin}. He first requires
every limit at infinity of $L_f$ to lie in the Steinberg variety
\[
 \ZSt=\{(-p,q)\in T^*(\Fl\times\Fl):\mu_\C(p)=\mu_\C(q)\}.
\]
A limit at infinity is the limit of a sequence of covectors after
division by their norms, as the norms tend to infinity.
By the moment-map argument in Lemma~\ref{lem:induction}, $f$ preserves
$[k,\xi]\mapsto\Ad_k\xi$ and hence the invariant norm
$|[k,\xi]|=|\xi|$. Homogeneity of
$\mu_\C$ and compactness of the unit cotangent bundle now identify
this graph condition with~\eqref{eq:end-rescaled}.

Jin further requires $L_f$ to lie in the diagonal conormal
$\{(-p,p):p\in T^*\Fl\}$ near each of its directions at infinity
belonging to no other component of $\ZSt$. These are exactly the
directions for which $\mu_\C(p)$ is regular. Indeed, write
$p=(V_\bullet,Z)$, where $ZV_j\subset V_{j-1}$. The components of
$\ZSt$ are the closures of the conormal bundles to the orbits of the
diagonal $SL_N(\C)$-action on pairs of flags. Over the open set of
regular nilpotents the flag is unique, with $V_j=\ker Z^j$, so only
the diagonal component occurs. If $Z$ is nonregular, some entry
$Z_{i,i+1}$ vanishes in a basis adapted to $V_\bullet$. Then
$ZV_{i+1}\subset V_{i-1}$, so one can vary $V_i$ between $V_{i-1}$
and $V_{i+1}$ while keeping $Z$ fixed. The resulting pairs of flags
show that $(-p,p)$ belongs to the conormal closure of the orbit of
pairs differing only at $i$.

By the first graph condition, every graph limit over a regular input
direction is diagonal, since its Springer fiber is a point. Thus if
exact fixation failed arbitrarily far out over a compact set of regular
unit covectors, off-diagonal graph points would accumulate at one of
these diagonal directions. Conversely, a cutoff on a compact
neighborhood of each regular unit covector gives Jin's condition
there. This identifies his second condition with that of
Definition~\ref{def:end}.

For convergence Jin requires compact-open $C^\infty$ convergence and
the first graph condition for limits along the varying graphs
\cite[(2.1)]{Jin}. On a compact parameter interval the latter is
exactly the uniform estimate~\eqref{eq:end-family}. Indeed, its
failure gives a sequence $s_n\to s$ and $|p_n|\to\infty$ for which a
normalized graph limit lies outside $\ZSt$. The cutoff in regular
directions is imposed separately on each map.

To compare twist conventions, use the trace normalization of
Lemma~\ref{lem:explicit-twist} in both coordinates. At the standard
flag, let $\xi^+$ be the strictly upper triangular part of the
skew-Hermitian representative $\xi$. The holomorphic cotangent
coordinate for the trace pairing is $u=-2\xi^+$, since
$\operatorname{Re}\tr(uA)=-\tr(\xi A)$ for skew-Hermitian $A$.
In these units Jin's Hermitian coordinate is
$\mathrm i(u-u^*)/2=-\mathrm i\xi$, where $u^*$ is the conjugate
transpose. His skew matrix on the $i$th root block is therefore
$-Q(z_i)$. Use the radial coordinate $r=|z_i|$ of~\eqref{eq:familytwist}
and take Jin's profile to be $h(r)=\pi-\vartheta(r)$ for $r\ge0$,
with $h(-r)=\pi-h(r)$. A positive rescaling of the radial coordinate
is absorbed by the corresponding rescaling of the argument of $h$.

For $r_i=|z_i|>0$ the multiplier in \cite[(3.5)]{Jin} is
\[
 \exp\!\left((\vartheta(r_i)-\pi)\frac{Q(z_i)}{r_i}\right)
 =\exp\!\left(\vartheta(r_i)\frac{Q(z_i)}{r_i}\right)(-I_2).
\]
The factor $-I_2$ on this block belongs to $T$, so this gives the
same point of $K\times_T\talg^\perp$ as~\eqref{eq:familytwist}.
Equality extends to $r_i=0$ by continuity. Thus Jin's $i$th twist
with this profile is $\tau_i$.

\begin{lemma}\label{lem:flag-rank-two-relations}
Let $i,j$ be distinct vertices of $\Gamma$. If $i\not\sim j$, then
$\tau_i$ and $\tau_j$ commute. If $i\sim j$, then $\tau_i\tau_j\tau_i$
and $\tau_j\tau_i\tau_j$ are joined by a continuous path in
$\Endg_K(T^*\Fl)$.
\end{lemma}
\begin{proof}
For $\ell=i,j$, let
$t_\ell=\operatorname{Ind}_{H_\ell}^{H_{ij}}t_\vartheta$ on
$T^*(H_{ij}/T)$. By Lemma~\ref{lem:induction},
$\tau_\ell=\operatorname{Ind}_{H_{ij}}^K t_\ell$, and induction
respects composition. If $i\not\sim j$, then $H_{ij}/T=S^2\times S^2$
by~\eqref{eq:small-flag-fibers}. Each $t_\ell$ twists its own
cotangent factor and fixes the other, so the maps commute.

If $i\sim j$, then $H_{ij}/T=SU(3)/T^2$, with $T^2$ the diagonal
maximal torus. Label its two simple roots by $i,j$ in that order.
Theorem~\ref{thm:jin-braid-map} gives
\[
 \beta_{SU(3)}(t_it_jt_i)=\sigma_1\sigma_2\sigma_1
 =\sigma_2\sigma_1\sigma_2=\beta_{SU(3)}(t_jt_it_j).
\]
Since $\beta_{SU(3)}$ is injective on path components, the two
products are joined by a path satisfying~\eqref{eq:end-family}.
The action of $H_{ij}$ on the flag fiber is generated by the image
of $SU(3)$ and a central torus acting trivially. Hence the path is
$H_{ij}$-equivariant.
Lemma~\ref{lem:induction-end} extends it to the required path on
$T^*\Fl$.
\end{proof}

\begin{proof}[Proof of Theorem~\ref{main:flag}]
The twists satisfy Definition~\ref{def:end} by
Lemma~\ref{lem:induction-end}, since $t_\vartheta$ is compactly
supported. Lemma~\ref{lem:flag-rank-two-relations} gives the defining
Artin relations. In type $A_{N-1}$ the cover $SU(N)\to K$ identifies
the flag cotangent bundles and equivariant groups.

For faithfulness we check that Jin's map is constant on paths.
For $N=2,3$ this follows from Theorem~\ref{thm:jin-braid-map}.
For $N\ge4$, let $f_s$, $0\le s\le1$, be a path and use the marked
disk in Jin's construction \cite[Lemma~3.2 and \S\,3.2]{Jin}.
In his Hermitian coordinates let $a=\operatorname{diag}(1,-1,0,\ldots,0)$
and let $M_a=\mu^{-1}(a)/K_a$, where $\mu$ is the compact moment map
and $K_a$ is the stabilizer of $a$. In a product chart in $M_a$, Jin's
disk contains $N$ marked points joined consecutively by arcs. Their
chain is an entire fiber of $[p]\mapsto[\mu_\C(p)]$ from $M_a$ to
$\mathcal N/K_a$. Two more arcs join the ends of the chain to distinct
points on the disk boundary. Each point of these outer arcs away from
the chain has singleton fiber. The relative homotopy classes of the
images of all these arcs, projected to the disk at a sufficiently
large moment level, determine the braid.

Identify $M_{ra}$ with $M_a$ by cotangent dilation. By
\eqref{eq:end-family}, every limit of pairs
$(p,r^{-1}f_s(rp))$ for $p\in\mu^{-1}(a)$ and $r\to\infty$ has
equal Springer images, even when $s$ varies. Since $M_a$ is compact
and the union of the arcs contains the entire fiber through each of
its points, the rescaled images of the arcs stay in a slightly
enlarged product chart for all sufficiently large $r$, uniformly in $s$.
At the two regular outer endpoints the same argument keeps the
projections in disjoint small disks containing no marks. Take $r$ large enough also to compute the braids of $f_0$ and
$f_1$ in Jin's construction. In the chart the marked strata are
the points whose stabilizers contain a circle. Equivariance keeps the
projected arc interiors away from the marks, and continuity keeps each
marked endpoint at the same mark. Join each moving outer endpoint to
its original position in its small disk and retract the enlarged
disk collar to the original boundary. This gives homotopies of the
projected arcs relative to their endpoints. Hence the two endpoint
maps give the same braid.

Write $\overline\beta_K$ for the induced map on path components.
By the identification of the twists in Theorem~\ref{thm:jin-braid-map},
\[
 Br_N\xrightarrow{\ \rho_{\rm flag}\ }
 \pi_0\Endg_K(T^*\Fl)\xrightarrow{\ \overline\beta_K\ }Br_N,
 \qquad \sigma_i\longmapsto\sigma_i .
\]
Thus $\rho_{\rm flag}$ is injective.
\end{proof}

\section{Entropy of twist words on tree plumbings}
\label{sec:tree}

In this section we prove Theorem~\ref{main:tree}. We count intersections
of curves with their iterates under $\phi_\beta$ on $\Sigma_\Gamma$.
Keating's comparison transfers this growth to Floer ranks in $W_\Gamma$.
The resulting area growth bounds the entropy of every representative
of $[\Phi_\beta]$ from below.

We first recall how to plumb two disk cotangent bundles.
\begin{definition}[{Cotangent plumbing, see e.g.\cite[Definition~2.1]{AbouzaidSmithPlumbings}}]
Let $Q_1,Q_2$ be closed smooth manifolds of dimension $m$. Choose a
coordinate disk in each $Q_i$ and identify the corresponding cotangent
blocks $D^m\times D^m\subset D^*Q_i$ by $(q,p)\mapsto(p,-q)$.
The resulting Liouville domain, after rounding corners, is the
\emph{plumbing} of $D^*Q_1$ and $D^*Q_2$. Its zero sections are exact
Lagrangians meeting transversely at one point.
\end{definition}

For a finite tree $\Gamma$, take one copy of $D^*S^2$ per vertex and
plumb the copies along the edges, using disjoint disks at each vertex.
The resulting domain is $W_\Gamma$, with \emph{core spheres} $L_v$
given by the zero sections.
Plumbing annuli $D^*S^1$ in the same way gives $\Sigma_\Gamma$, whose
zero sections are the \emph{core curves} $\gamma_v$. In either plumbing
two cores meet once if their vertices are adjacent and are otherwise
disjoint \cite[Example~1.2]{KeatingStab}.
For the annular model on a chain see \cite[\S\,3.3]{GWX}.

Both unmarked domains depend only on $\Gamma$ up to Liouville
deformation. Indeed, the surface has Euler characteristic $1-|V(\Gamma)|$
and twice its genus is the rank of the core intersection matrix, whose
edge signs can be fixed recursively on the tree. To specify the marked
curves $\gamma_v$, we fix the cyclic order of the disks on each annulus
\cite[Example~1.2]{KeatingStab}.

In both plumbings the Dehn twists about disjoint zero sections commute,
and twists about zero sections meeting once satisfy the braid relation
up to symplectic isotopy \cite{SeidelBook,KeatingStab}. Thus the words in
\eqref{eq:arbitrary-tree-word} define mapping classes $\phi_\beta$ and
$[\Phi_\beta]$ for each $\beta\in\Art_\Gamma$. The twists are supported
in the interior. Each isotopy has support in a common compact subset
of the interior, and we measure entropy on $W_\Gamma$.

\subsection{Growth of intersections on the surface}

\begin{theorem}[{Nielsen--Thurston classification
\cite{ThurstonSurfaces}, \cite[Expos\'e~11]{FLP}}]
Let $h$ represent an orientation-preserving mapping class of a connected
oriented surface of finite type and negative Euler characteristic, with
boundary circles treated as punctures. The class is periodic,
pseudo-Anosov, or reducible.
In the reducible case $h$ can be chosen to preserve the canonical
multicurve of disjoint essential nonperipheral simple closed curves.
After cutting along it
and replacing the cut circles by punctures, some power of $h$ fixes
each component and is periodic or pseudo-Anosov on it.
\end{theorem}

A periodic class has finite order. A pseudo-Anosov representative
scales two transverse measured foliations by $\lambda$ and
$\lambda^{-1}$, where $\lambda>1$ is its stretch factor.

\begin{definition}\label{def:pA-growth}
For the annulus set $\Lambda_\beta=1$. Otherwise let $h$ be a
Nielsen--Thurston representative of $\phi_\beta$. For each component
$Y$ of its decomposition, let $p_Y=\min\{p\ge1:h^p(Y)=Y\}$.
If the first return $h^{p_Y}|_Y$ is pseudo-Anosov, denote
its stretch factor by $\lambda_Y$. Set
\begin{equation}\label{eq:pA-growth}
 \Lambda_\beta=\max\bigl(\{1\}\cup
 \{\lambda_Y^{1/p_Y}:Y\text{ has pseudo-Anosov first return}\}\bigr).
\end{equation}
\end{definition}

Write $i(a,b)$ for the geometric intersection number of $a$ and $b$.
\begin{theorem}[{Thurston's intersection growth formula
\cite[Expos\'e~12, Theorem~12.2]{FLP}}]
\label{thm:tree-surface-growth-input}
If $\psi$ is pseudo-Anosov on a closed orientable surface of genus at least two with stretch factor $\lambda$, then for any
essential nonperipheral simple closed curves $a,b$ on its surface,
\[
 \lim_{n\to\infty} i(a,\psi^n b)^{1/n}=\lambda .
\]
\end{theorem}

\begin{lemma}\label{input:surface-growth}
Let $\gamma_i$ be the zero section of the $i$th annulus in the
plumbing construction of $\Sigma_\Gamma$. Then
\begin{equation}\label{eq:surface-geometric-growth}
 \limsup_{n\to\infty}\frac1n\log\left(1+
 \sum_{i,j}i(\gamma_i,\phi_\beta^n\gamma_j)\right)
 \ =\ \log\Lambda_\beta .
\end{equation}
Moreover $\rho(M_\beta)\le\Lambda_\beta$.
\end{lemma}
For more than one vertex this is Thurston's spectral theorem applied
to the filling family of cores \cite[Expos\'e~11, Theorem~11.9 and its
proof]{FLP}. Total intersection with the cores is comparable to geodesic
length, which bounds the homology norm in the core basis. For the
annulus the intersection sum is zero and $M_\beta=I$.

\subsection{Floer ranks in the plumbing}

We use Keating's Lefschetz stabilization to compare the sphere twists
in $W_\Gamma$ with the corresponding surface twists in $\Sigma_\Gamma$.

\begin{lemma}\label{prop:tree-floer-comparison}
For all vertices $i,j$ and every $n\ge0$,
\[
 \dim_{\mathbb F_2}\HF_{W_\Gamma}(L_i,\Phi_\beta^nL_j)
                 \ \ge\ i(\gamma_i,\phi_\beta^n\gamma_j).
\]
\end{lemma}
\begin{proof}
In Keating's one-stabilization $M$ of $\Sigma_\Gamma$, each core occurs
twice in the ordered list of vanishing cycles. Pair its two copies by a
diameter of the base disk \cite[Definition~1.3]{KeatingStab}. The matching
spheres meet only at core crossings in the central fiber. There the
core directions and the two diameter directions give transversality.
A cotangent neighborhood of this configuration is the plumbing $W_\Gamma$.
Its primitive and the ambient primitive differ by $dh$, since
$H^1(W_\Gamma;\R)=0$. Subtracting the differential of a cutoff extension
of $h$ makes the primitives agree near $\partial W_\Gamma$.

Exact Floer cohomology of compact exact Lagrangians in this subdomain
agrees with that in $M$ \cite[proof of Proposition~2.9]{KeatingStab}.
Keating's equivariant twists restrict to the surface twists
\cite[Lemma~2.6]{KeatingStab}. They are symplectically isotopic to the
corresponding plumbing twists extended by the identity. Thus their
word $\widetilde\Phi_\beta$ and the extended word $\Phi_\beta$ have
Hamiltonian isotopic images of each sphere, since $H^1(S^2;\R)=0$.
Keating's comparison \cite[Theorem~2.7]{KeatingStab} for the inverse
word gives, over $\mathbb F_2$,
\begin{align*}
 \dim\HF_{W_\Gamma}(L_i,\Phi_\beta^nL_j)
 &=\dim\HF_M(\widetilde\Phi_\beta^{-n}L_i,L_j)\\
 &\ge\dim\HF_{\Sigma_\Gamma}(\phi_\beta^{-n}\gamma_i,\gamma_j)
 =\dim\HF_{\Sigma_\Gamma}(\gamma_i,\phi_\beta^n\gamma_j).
\end{align*}
For nonisotopic curves the last rank equals
$i(\gamma_i,\phi_\beta^n\gamma_j)$ by the surface Floer computation of
de Silva, Robbin and Salamon
\cite[Theorem~9.3 and Corollary~9.5]{deSilvaRobbinSalamon}, applied
on $\Int\Sigma_\Gamma$. For isotopic curves the intersection number
is zero.
\end{proof}

\subsection{From Floer ranks to entropy}\label{subsec:floer-entropy}

In the Floer-theoretic form of Arnol'd's intersection estimate, Floer
rank bounds transverse intersection counts from below
\cite[Lemma~5.14 and the discussion preceding Proposition~5.15]{SeidelLecture}.
The Crofton inequality turns the average of these counts over Hamiltonian perturbations
into a lower bound on the area of each iterated sphere. By Yomdin's theorem,
the exponential growth rate of this area bounds entropy from below.

\begin{theorem}[{\cite[Lemmas~5.3 and~5.6]{CGG}}]
\label{thm:tree-crofton}
Let $P$ be a closed Lagrangian in the interior of a symplectic
$2m$-manifold $M$. If $P$ admits an immersion into $\R^d$, then
there is a small $d$-parameter family of Hamiltonian images
$P_s$ of $P$, parametrized by a compact ball $B$, with evaluation
$B\times P\to M$ extending to a submersion on a neighborhood of
$B\times P$ in $\R^d\times P$.
The perturbations can be supported in a cotangent neighborhood of $P$.
For any smooth density $ds$ on $B$ and any compact smooth
$m$-submanifold $S$, almost every $P_s$ meets $S$ transversely and
\[
 \int_B\#(P_s\cap S)\,ds\le C\operatorname{vol}_m(S),
\]
where $C$ depends on the family, density and ambient metric.
\end{theorem}

\begin{theorem}[{Yomdin's volume bound \cite{Yomdin}}]
\label{thm:tree-yomdin}
Let $M$ be a closed smooth Riemannian manifold and
$f:M\to M$ a $C^\infty$ diffeomorphism. For each compact smooth
embedded submanifold $S$ of dimension $m$, possibly with boundary,
\[
 \limsup_{n\to\infty}\frac1n
     \log\operatorname{vol}_m(f^n S)\le h_{\rm top}(f).
\]
The disk version implies this form by covering $S$ by finitely many
smooth parameterized disks and adding their volume bounds.
\end{theorem}

\begin{proposition}\label{prop:floer-entropy}
Let $M$ be a compact simply connected symplectic four-manifold,
possibly with boundary.
Fix finite collections of compact Lagrangian spheres
$P_1,\ldots,P_a$ and $Q_1,\ldots,Q_b$ in its interior.
Assume the relevant Floer groups are finite-dimensional,
Hamiltonian invariant, and bounded in rank by transverse
intersections. For a smooth symplectomorphism $g$, fixed near the
boundary when present, put
\[
 R_n(g)=\sum_{i=1}^a\sum_{j=1}^b
                       \dim\HF_M(P_i,g^nQ_j).
\]
If $f$ represents the same symplectic class, with the common-support
convention in the boundary case, then $R_n(f)=R_n(g)$ and
\begin{equation}\label{eq:floer-entropy}
 \htop(f)\ge\limsup_{n\to\infty}\frac1n\log\bigl(1+R_n(g)\bigr).
\end{equation}
\end{proposition}

\begin{proof}
For each $n$, choose a symplectic isotopy from $g^n$ to $f^n$.
Its generating one-form is exact since $H^1(M;\R)=0$. If $M$ has
boundary, multiply the generating Hamiltonian by a cutoff equal to
$1$ near the compact trace of the spheres $Q_j$ and supported in
$\Int M$. The images $g^nQ_j$ and $f^nQ_j$ are thus Hamiltonian
isotopic, so $R_n(f)=R_n(g)$.

Each $P_i$ is a sphere and embeds in $\R^3$. Use the Hamiltonian
family $P_{i,s}$ of Theorem~\ref{thm:tree-crofton}, supported in a
cotangent neighborhood, with parameter measure of mass one. It is
transverse to $f^nQ_j$ for almost every $s$. The Floer rank bound
followed by Crofton gives
\[
 R_n(g)\le\sum_{i,j}\int\#(P_{i,s}\cap f^nQ_j)\,ds
              \le\Big(\sum_i C_i\Big)\sum_j\Area(f^nQ_j),
\]
where the $C_i$ are independent of $n$. If $M$ has boundary, extend
$f$ by the identity on the other half of the double. The iterated
spheres and the entropy are unchanged. Theorem~\ref{thm:tree-yomdin} gives
\[
 \limsup_{n\to\infty}\frac1n\log\Area(f^nQ_j)\le h_{\rm top}(f)
 \qquad(1\le j\le b).
\]
The sum of these areas and the added constant $1$ obey the same
upper growth bound, proving \eqref{eq:floer-entropy}.
\end{proof}

\begin{proof}[Proof of Theorem~\ref{main:tree}]
The plumbing retracts onto a tree of spheres and is therefore simply connected. Exact Floer cohomology over $\mathbb F_2$
is Hamiltonian invariant and satisfies the intersection bound.
Apply Proposition~\ref{prop:floer-entropy} with both collections
equal to the cores. Lemmas~\ref{prop:tree-floer-comparison}
and~\ref{input:surface-growth} give
\[
 h_{\rm top}(f)\ge
 \limsup_{n\to\infty}\frac1n\log\left(1+\sum_{i,j}
        i(\gamma_i,\phi_\beta^n\gamma_j)\right)
 =\log\Lambda_\beta\ge\log\rho(M_\beta).
\]
\end{proof}

\subsection{Penner words and explicit entropy bounds}
\label{subsec:penner-tree}

For a word using every vertex, with opposite twist signs on the two parts
of the tree, Penner's construction computes the entropy bound
$\log\Lambda_\beta$ from the intersection matrix of the core curves.

\begin{theorem}[{Penner's construction \cite[Theorem~3.1]{Penner},
as used in \cite[Theorem~2.4]{GWX}}]
\label{thm:penner}
Let $S$ be a connected orientable surface of finite type with
negative Euler characteristic, and allow punctures. Let $A$ and $B$
be multicurves in minimal position that fill $S$: the complementary
regions of $A\cup B$ are disks or once-punctured disks.
A product of positive Dehn twists about curves of $A$ and negative
Dehn twists about curves of $B$, using every curve at least once,
represents a pseudo-Anosov mapping class.

Order the curves of $A\cup B$ as $c_1,\ldots,c_d$. Put
$\Omega_{ij}=i(c_i,c_j)$ and let $D_i$ be the diagonal matrix with
entry $1$ at $(i,i)$ and zero elsewhere. The stretch factor is the
Perron--Frobenius eigenvalue of the corresponding product of
$Q_i=I+D_i\Omega$, with one factor $Q_i$ for each occurrence of the
allowed twist about $c_i$, in the same composition order
\cite[\S\,2.1]{ShinStrenner}.
\end{theorem}

\begin{corollary}
\label{cor:penner-tree}
Suppose $\Gamma$ is a finite tree with at least two vertices and
bipartition $A\sqcup B$. Let $\beta$ be a word using every vertex,
with positive powers at vertices of $A$ and negative powers at
vertices of $B$. Then
\[
 \Lambda_\beta=\rho(M_\beta)>1.
\]
In particular, every representative $f$ of $[\Phi_\beta]$ in
$\Symp_c(\Int W_\Gamma)$ satisfies $\htop(f)\ge\log\rho(M_\beta)$.

Let $C$ be the adjacency block from $A$ to $B$, so that
$C_{ab}=1$ on the edges and $0$ otherwise, and set
$\sigma=\|C\|_2$. For the two-block word with surface representative
\[
 \phi_\beta=\Bigl(\prod_{a\in A}T_{\gamma_a}^{\,r}\Bigr)
             \Bigl(\prod_{b\in B}T_{\gamma_b}^{-s}\Bigr),
 \qquad r,s\in\Z_{>0},
\]
the stretch factor is the larger root of
\begin{equation}\label{eq:penner-two-block}
 \Lambda_\beta+\Lambda_\beta^{-1}=2+rs\sigma^2 .
\end{equation}
\end{corollary}

\begin{proof}
The two colors give disjoint multicurves within each color.
The cores fill the punctured surface and are in minimal position,
since neighboring cores cross only once. Penner's theorem
(Theorem~\ref{thm:penner}) gives a pseudo-Anosov class on this surface,
whose stretch factor is $\Lambda_\beta$.

Orient the cores recursively along edges so that, for $a\in A$ and
$b\in B$ adjacent, the two twists act by
\[
 T_{\gamma_a}[\gamma_b]=[\gamma_b]+[\gamma_a],
 \qquad T_{\gamma_b}^{-1}[\gamma_a]=[\gamma_a]+[\gamma_b].
\]
The second formula follows from the first by skew symmetry of
algebraic intersection and the sign of the inverse twist. The
matrix of each allowed twist in this basis is $Q_v$ from
Theorem~\ref{thm:penner}. Thus $M_\beta$ is the Penner product,
$\rho(M_\beta)=\Lambda_\beta>1$, and Theorem~\ref{main:tree}
gives the entropy bound.

For the two-block word the product is
\[
 M_\beta=
 \begin{pmatrix}I&rC\\0&I\end{pmatrix}
 \begin{pmatrix}I&0\\sC^{\mathsf T}&I\end{pmatrix}
 =
 \begin{pmatrix}I+rsCC^{\mathsf T}&rC\\sC^{\mathsf T}&I\end{pmatrix}.
\]
The adjacency matrix
$\left(\begin{smallmatrix}0&C\\C^{\mathsf T}&0\end{smallmatrix}\right)$
has a positive Perron eigenvector $(u,v)$ with
$Cv=\sigma u$ and $C^{\mathsf T}u=\sigma v$. On the basis
$(u,0),(0,v)$ the restriction of $M_\beta$ is
\[
 \begin{pmatrix}1+rs\sigma^2&r\sigma\\s\sigma&1\end{pmatrix}.
\]
Its larger eigenvalue has a positive eigenvector in the full core
basis, so it equals $\rho(M_\beta)$. The determinant is $1$ and
the trace is $2+rs\sigma^2$, proving \eqref{eq:penner-two-block}.
\end{proof}

\section{Slodowy fibers and transverse monodromy}
\label{sec:slodowy}

In this section we prove Theorem~\ref{main:slice}. We identify the
monodromy of an Artin word with the corresponding product of Dehn
twists about the core spheres of the ADE plumbing. The entropy bound
then follows directly from Theorem~\ref{main:tree}. Faithfulness
follows from Khovanov--Seidel in type $A$ and from Brav--Thomas's
theorem on spherical twists in types $D$ and $E$.

\subsection{Meridians and vanishing cycles}\label{subsec:flag-slice-comparison}

Recall that $F_\Gamma$ is the part of the regular fiber $X_{b_*}$
that lies in the Slodowy slice $\mathcal S_e$ and inside a small ball
about $e$, as in \eqref{eq:compact-slice-fibration}. We follow the
same loops as in Theorem~\ref{main:flag}, taking parallel transport
in the slice and fixing a collar of the boundary.

The completions of $F_\Gamma$ and the plumbing $W_\Gamma$ are
symplectomorphic. In this model the sphere that collapses along the
vanishing path for vertex $v$ is the core sphere $L_v$, the zero
section of the $v$th disk cotangent bundle $D^*S^2$.%
\footnote{For the plumbing model see Y.~Lekili and K.~Ueda,
\emph{Homological mirror symmetry for Milnor fibers of simple
singularities}, Algebr. Geom. \textbf{8} (2021), 562--586, \S\,3.
For the vanishing spheres see \cite[Lemma~18.2]{SeidelBook} and
\cite[Lemma~15.9]{SeidelLecture}.}

\begin{theorem}[{Brieskorn--Slodowy
\cite{Brieskorn}, \cite[\S\S\,6.2 and~8.7]{Slodowy}}]
Let $e$ be subregular nilpotent in a complex simple Lie algebra of
type $\Gamma$ in ADE. The germ of $\chi_{\Gamma,e}^{-1}(0)$ at $e$
is the simple surface singularity of type $\Gamma$, and
$\chi_{\Gamma,e}$ is its semiuniversal deformation. Its discriminant
is the image of the reflection hyperplanes in $\halg/\mathsf W_\Gamma$.
At a smooth discriminant point the fiber has one ordinary double
point. Over a transverse disk the local model at that point is
$z_1^2+z_2^2+z_3^2=t$.
\end{theorem}

\begin{theorem}[{Symplectic Picard--Lefschetz monodromy
\cite[Proposition~1.15]{SeidelLES}, \cite{SeidelBook}}]
\label{thm:slice-picard-lefschetz}
The monodromy around a positive meridian at a smooth point of the
discriminant is the positive Dehn twist in the associated vanishing
sphere $L$, up to symplectic isotopy fixed near the fiber boundary.
On $H_2(F_\Gamma)$ it acts by $D\mapsto D+(D\cdot[L])[L]$,
where $[L]^2=-2$.
\end{theorem}

The following standard description of slice monodromy combines
Brieskorn's presentation \cite{BrieskornBraid} with the boundary
normalization of Khovanov--Seidel
\cite[Proposition~6.2 and Lemma~6.8]{KhovanovSeidel}. We include
a proof sketch for completeness.

\begin{lemma}\label{lem:slice-monodromy}
The boundary-normalized monodromy of \eqref{eq:compact-slice-fibration}
defines a homomorphism
\[
 \rho_{\mathrm{slice}}\colon\Art_\Gamma\longrightarrow
 \pi_0\Symp_c(\Int F_\Gamma),\qquad \sigma_v\longmapsto[\tau_{L_v}].
\]
\end{lemma}
\begin{proof}[Proof sketch]
The reflection hyperplanes are linear, so radial rescaling identifies
the fundamental groups of the local and global discriminant complements.
Brieskorn's presentation identifies the simple meridians with the
generators $\sigma_v$ of $\Art_\Gamma$.\footnote{For the distinguished
paths see K.~Saito, \emph{Polyhedra dual to the Weyl chamber
decomposition: a pr\'ecis}, Publ. Res. Inst. Math. Sci. \textbf{40}
(2004), 1337--1384, \S\,4.}
The fiber boundaries form a smooth contact family over the small
contractible neighborhood $\mathcal V_\Gamma$, including over the
discriminant. Choose a contact trivialization of this family and deform the Liouville forms of the regular fibers on a collar, keeping the reference fiber fixed, so that the trivialization identifies the boundary contact forms. The resulting strict symplectic fibration admits product transport on a collar~\cite[Lemmas~6.3, 6.4 and~6.8]{KhovanovSeidel}. Thus based loops
give compactly supported mapping classes, and
Theorem~\ref{thm:slice-picard-lefschetz} sends the $v$th meridian
to $[\tau_{L_v}]$.
\end{proof}

\subsection{Faithfulness in type \texorpdfstring{$A$}{A}}

In type $A_r$ the slice has the polynomial model below.
Its simple meridians are the half-twists exchanging adjacent roots,
so the slice action is the braid action studied by Khovanov and Seidel
\cite[\S\,6b]{KhovanovSeidel}.

\begin{theorem}[{Khovanov--Seidel
\cite[Corollary~1.4 and \S\,6c]{KhovanovSeidel}}]
\label{thm:ks-slice-faithful}
Let $r\ge1$. For the surface Milnor family $x^2+y^2+P(z)=0$, where
$P$ is a small monic deformation of $z^{r+1}$ with simple roots, the
boundary-normalized monodromy
$Br_{r+1}\to\pi_0\Symp_c(\Int F_{A_r})$ is injective. A positive
half-twist exchanging adjacent roots acts by the Dehn twist in the
matching sphere over the joining arc.
\end{theorem}

We postpone the proof of faithfulness in types $D$ and $E$ to
Appendix~\ref{app:slice-faithful-de}, where we apply Brav--Thomas's
theorem to the graded vanishing spheres \cite[Theorem~3.1]{BravThomas}.

\begin{proof}[Proof of Theorem~\ref{main:slice}]
By Lemma~\ref{lem:slice-monodromy}, the monodromy of $\beta$ is the
class of the twist word $\Phi_\beta$ in \eqref{eq:arbitrary-tree-word}.
Theorem~\ref{thm:ks-slice-faithful} and Lemma~\ref{lem:slice-faithful}
give faithfulness.

Let $f$ represent $\rho_{\mathrm{slice}}(\beta)$. Extend $f$ by the
identity to the completion and use the plumbing identification above.
Choose a compact truncation containing the supports of $f$ and its
isotopy to the core twist word. Liouville rescaling carries these
supports into $W_\Gamma$ and preserves the twist mapping classes
\cite[Remark~6.5]{KhovanovSeidel}. Conjugation and extension by the
identity preserve entropy, so Theorem~\ref{main:tree} gives
\[
 h_{\rm top}(f)\ge\log\Lambda_\beta\ge\log\rho(M_\beta).
\]
\end{proof}

For a Penner word as in Corollary~\ref{cor:penner-tree}, the surface
class is pseudo-Anosov and every representative of
$\rho_{\mathrm{slice}}(\beta)$ has entropy at least
$\log\rho(M_\beta)>0$. For the two-block word with exponents $r,-s$,
let $\sigma$ be the largest singular value of the adjacency block.
The larger root in \eqref{eq:penner-two-block} gives
\[
 h_{\rm top}(f)\ge
 \log\!\left(1+\frac{rs\sigma^2}{2}
             +\frac{\sqrt{rs\sigma^2(4+rs\sigma^2)}}{2}\right).
\]
For $r=2p$ and $s=2q$ this is $\log\lambda_\Gamma(p,q)$
in \eqref{eq:K3-lambda}.

\section{Global Floer growth on K3 surfaces}\label{sec:k3-main}
\label{sec:floer}

We prove Theorem~\ref{main:k3} and Corollary~\ref{cor:k3-A2} by
computing Floer ranks in the closed K3 surface $X$. The twists are
supported in a plumbing neighborhood of the ADE configuration. Their
action on the Floer ranks of a suitable sphere gives the growth rate.

In this section, we use the universal Novikov field
\begin{equation}\label{eq:Novikov}
 \Lambda=\left\{\sum_{\alpha\in\R}c_\alpha T^\alpha:
 c_\alpha\in\C,\quad
 \#\{\alpha<E:c_\alpha\ne0\}<\infty\text{ for every }E\in\R\right\},
 \qquad |T|=0 ,
\end{equation}
and denote by $\HF_X(L_1,L_2)$ the Floer cohomology of the Lagrangian
spheres $L_1,L_2\subset X$ over $\Lambda$.
Let $A_1,\ldots,A_m$ and $B_1,\ldots,B_\ell$ be the two colors of the
configuration in Theorem~\ref{main:k3}, with adjacency block $C$.
For a Lagrangian sphere $K\subset X$, write its total Floer ranks as
\[
 x(K)=\bigl(\dim_\Lambda\HF_X(A_i,K)\bigr)_{i=1}^m,\qquad
 y(K)=\bigl(\dim_\Lambda\HF_X(B_j,K)\bigr)_{j=1}^{\ell}.
\]
A \emph{twist descendant} is a sphere obtained from a core sphere by
successive positive or negative Dehn twists about the core spheres.

\begin{proposition}\label{prop:rank-recurrence}
Choose core spheres $A_a$ and $B_b$ that meet once. Let
$K_*=\tau_{A_a}^4(B_b)$ be the sphere obtained by twisting $B_b$ four
times about $A_a$. For every integer $n\ge0$, the Floer ranks satisfy
\begin{equation}\label{eq:K3-matrix}
 \binom{x(\Phi_{p,q}^{\,n}K_*)}{y(\Phi_{p,q}^{\,n}K_*)}
 =
 \begin{pmatrix} I_m&2qC\\2pC^t&I_\ell+4pqC^tC\end{pmatrix}^{\!n}
 \binom{x(K_*)}{y(K_*)} .
\end{equation}
\end{proposition}

In the next few subsections, we shall prove Proposition~\ref{prop:rank-recurrence}
using Seidel's exact triangle. Starting from $K_*$, we show that the Floer
product along each incident edge vanishes at every twist. The triangle
then gives the matrix in \eqref{eq:K3-matrix}.

\subsection{The global twist sequence}\label{subsec:comparison-reduction}

For transverse spheres $F,K$, the Floer complex has one
generator per intersection point, so $\dim\HF_X(F,K)\le\#(F\cap K)$.

Fix a trivialization of the squared canonical bundle of $X$. Each
Lagrangian sphere admits a grading since $H^1(S^2;\Z)=0$. Capping a
disc along its boundary in the sphere identifies its Maslov index
with twice the Chern number of the capped sphere, which is zero on $X$.
We use the sphere's
unique spin structure and trivial local system. Write
$[\omega]=c\,c_1(\mathcal H)$ as in Definition~\ref{def:polarized-k3}.
A unitary connection on $\mathcal H$ with curvature
$-2\pi\mathrm i\,\omega/c$ restricts to a flat connection with
trivial holonomy on each sphere. Thus the spheres
are rational branes in Seidel's closed Fukaya category. We use his
regular boundary data and the identification of graded isotopic
spheres \cite[Lemmas~8.4 and~8.14]{SeidelQuartic}. A small Hamiltonian
pushoff of $L$ using a Morse function with two critical points gives
one generator in each of degrees zero and two. The differential is
therefore zero, and
\begin{equation}\label{eq:self}
 \HF_X^*(L,L)=\Lambda 1_L\oplus\Lambda u_L,\qquad |1_L|=0,\quad |u_L|=2 .
\end{equation}
Appendix~\ref{app:seidel-comparison} identifies the area weights for
$\omega/c$ with those for $\omega$ after extending scalars by
$q^a\mapsto T^{ca}$ for $a\in\Q$. We write the Floer triangle product as
composition $y\circ x$.

With $L[t]$ denoting a grading shift and the twist lift normalized
outside its support as in \cite[Lemma~5.7]{SeidelGraded}, we have
\begin{equation}\label{eq:shifts}
 \HF_X^k(L[t],K)=\HF_X^{k-t}(L,K),\qquad
 \HF_X^k(L,K[t])=\HF_X^{k+t}(L,K),\qquad \tau_LL=L[-1].
\end{equation}
A graded symplectomorphism $g$ transports the branes and preserves
the Floer product. Indeed, the pullback connection has the same
curvature, so the induced connection on
$\operatorname{Hom}(\mathcal H,g^*\mathcal H)$ is flat and has
trivial holonomy on the simply connected $X$. Applying naturality
to the two twists gives
\begin{equation}\label{eq:twist-shifts}
 \HF_X^k(A,\tau_AK)\cong\HF_X^{k-1}(A,K),\qquad
 \HF_X^k(B,\tau_B^{-1}K)\cong\HF_X^{k+1}(B,K).
\end{equation}

\begin{proposition}[{Seidel's closed twist sequence
\cite[(2.14)--(2.15) and Proposition~9.4]{SeidelQuartic}}]
\label{prop:closed-sequence}
For graded Lagrangian spheres $F,L,K$ in $X$, there is an exact sequence
\begin{equation}\label{eq:Seidel}
\begin{aligned}
 \cdots\longrightarrow(\HF_X(L,K)\otimes_\Lambda\HF_X(F,L))^k
 &\xrightarrow{\ \mu\ }\HF_X^k(F,K)\longrightarrow\HF_X^k(F,\tau_LK)\\
 &\longrightarrow(\HF_X(L,K)\otimes_\Lambda\HF_X(F,L))^{k+1}
 \longrightarrow\cdots ,
\end{aligned}
\end{equation}
where $\mu(y\otimes x)=y\circ x$ and the tensor product has its total grading.
\end{proposition}
Appendix~\ref{app:seidel-comparison} obtains \eqref{eq:Seidel} by
applying $\Hom^*(F,-)$ to the evaluation triangle
\[
 \HF_X^*(L,K)\otimes L\xrightarrow{\mathrm{ev}} K
 \longrightarrow\tau_LK\longrightarrow(\HF_X^*(L,K)\otimes L)[1].
\]
The second arrow induces maps on $\HF_X(F,-)$ that commute with
precomposition. Replacing $K$ by $\tau_L^{-1}K$ gives an arrow
$\tau_L^{-1}K\to K$ with the same property.
Writing $Q,G,D$ for the tensor product, $\HF_X(F,K)$ and
$\HF_X(F,\tau_LK)$ in \eqref{eq:Seidel}, exactness gives
\begin{equation}\label{eq:rankLES}
 \dim D=\dim Q+\dim G-2\rk\mu .
\end{equation}
Indeed, $D^k$ is an extension of $\ker\mu^{k+1}$ by
$\operatorname{coker}\mu^k$. If $\mu=0$, then
$\dim D^k=\dim G^k+\dim Q^{k+1}$.

The edge calculation uses the pairing given by composition and trace.
Ganatra gives a marked-disc description in the exact setting
\cite[Lemma~15]{GanatraCY}. For spheres in the closed K3 surface the
required pairing follows from the evaluation triangle.
\begin{lemma}\label{lem:k3-sphere-pairing}
Let $\tr_L\colon\HF_X^2(L,L)\to\Lambda$ send $u_L$ to $1$. For graded
Lagrangian spheres $L,K$ in $X$, the pairing
\begin{equation}\label{eq:pairing}
 \HF_X^k(L,K)\otimes_\Lambda\HF_X^{2-k}(K,L)\longrightarrow\Lambda,
 \qquad x\otimes y\longmapsto\tr_L(y\circ x)
\end{equation}
is perfect.
\end{lemma}
\begin{proof}
Put $V=\HF_X^*(L,K)$ and apply $\Hom^*(-,L)$
to the triangle. The resulting composition map is
\[
 \HF_X^*(K,L)\longrightarrow
 \Hom_\Lambda^*(V,\HF_X^*(L,L)),\qquad y\longmapsto(x\mapsto y\circ x).
\]
By \eqref{eq:shifts} and naturality, $\HF_X^*(\tau_LK,L)$ has
the same total dimension as
$\HF_X^*(K,L)$. The target of the composition map has dimension
$2\dim V$, so the rank identity for this triangle gives
$\dim V$ as the rank of the map.

If homogeneous classes satisfy $y\circ x=c1_L$ for $c\ne0$, then
$x\circ y$ is a scalar multiple of $1_K$ by \eqref{eq:self}. Composing
with $x$ shows that the scalar is $c$. In that case
$K$ is isomorphic to a shift of $L$ and \eqref{eq:pairing} is the
perfect pairing of the unit and top class in \eqref{eq:self}. In all
other cases the composition map takes values in
$\Hom_\Lambda^*(V,\Lambda u_L)$ and hence surjects onto this
$\dim V$-dimensional space. Interchanging $K$ and $L$ gives the
opposite dimension inequality, so the map is an isomorphism in
complementary degrees.
\end{proof}

\subsection{Products along one edge}\label{sec:rankproof}

Root the tree and grade each sphere so that the intersection generator
from its parent has degree one. By duality the reverse generator also
has degree one. We suppress the ambient subscript $X$. For an edge
$A$--$B$, choose generators $a\in\HF^1(A,B)$ and $b\in\HF^1(B,A)$.
Precomposition gives the degree-one maps
\begin{equation}\label{eq:alphabeta}
 \begin{aligned}
 \alpha_K:\HF(A,K)&\longrightarrow\HF(B,K),&x&\longmapsto x\circ b,\\
 \beta_K:\HF(B,K)&\longrightarrow\HF(A,K),&y&\longmapsto y\circ a.
 \end{aligned}
\end{equation}
For the edge $A_i$--$B_j$ we write $\alpha_{ij,K}$. Its vanishing is
independent of the generator. Write $R_u$ for precomposition with $u$.
By the perfect pairing \eqref{eq:pairing}, $b\circ a$ and $a\circ b$
are nonzero multiples of $u_A$ and $u_B$, respectively. Hence, up to
nonzero scalars,
\begin{equation}\label{eq:composites}
 \alpha_K\beta_K=R_{u_B},\qquad \beta_K\alpha_K=R_{u_A}.
\end{equation}

\begin{lemma}\label{lem:top}
If $R_{u_S}$ acts nontrivially on $\HF(S,K)$ for two graded spheres
$S,K$, then $K$ is isomorphic to a grading shift of $S$ and
$\dim\HF(S,K)=2$. If $S$ is a core sphere and $K$ a twist descendant,
then $[K]=\pm[S]$ in $H_2(X;\Z)$.
\end{lemma}
\begin{proof}
Take $x\in\HF^d(S,K)$ with $x\circ u_S\ne0$. By \eqref{eq:pairing}
there is $y\in\HF^{-d}(K,S)$ with
$\tr_S((y\circ x)\circ u_S)\ne0$. Thus $y\circ x$ is a nonzero
multiple of $1_S$. As in the proof of Lemma~\ref{lem:k3-sphere-pairing},
$x$ gives an isomorphism up to shift, and \eqref{eq:self} gives the rank.

For a twist descendant, Picard--Lefschetz puts $[K]$ in the span of
the core classes. The grading shift multiplies the Floer Euler
characteristics against all cores by a common sign $\varepsilon$, so
$[K]-\varepsilon[S]$ is orthogonal to every core. The intersection
form on their span is nondegenerate, since orienting the cores along
the tree gives the negative ADE Cartan matrix. Hence $[K]=\varepsilon[S]$.
\end{proof}

\begin{lemma}\label{lem:edge}
Suppose $\alpha_K=0$, and put $x=\dim\HF(A,K)$, $y=\dim\HF(B,K)$.
Twisting $K$ by $\tau_A$ changes these ranks to $(x,x+y)$, while
twisting by $\tau_B^{-1}$ changes them to $(x+y,y)$. At the new sphere
$K'$, $\alpha_{K'}=0$ if $R_{u_B}=0$ on $\HF(B,K')$ in the first
case, or $R_{u_A}=0$ on $\HF(A,K')$ in the second. In particular,
$x+y>2$ suffices in either case.
\end{lemma}
\begin{proof}
For any $K$, the rank identity \eqref{eq:rankLES} gives
\[
 \dim\HF(B,\tau_AK)=\dim\HF(A,\tau_B^{-1}K)
                  =x+y-2\rk\alpha_K.
\]
For $\tau_AK$, use the triple $(B,A,K)$ in \eqref{eq:Seidel}.
Its product is $\alpha_K$. For $\tau_B^{-1}K$, use $(K,B,A)$. Its product
$z\mapsto b\circ z$ is adjoint to $\alpha_K$ by \eqref{eq:pairing},
since $(x'\circ b)\circ z=x'\circ(b\circ z)$ for $x'\in\HF(A,K)$
and $z\in\HF(K,B)$. The third group in the triangle has rank
$\dim\HF(K,\tau_BA)=\dim\HF(A,\tau_B^{-1}K)$ by naturality and
duality. When $\alpha_K=0$, the displayed identity and
\eqref{eq:twist-shifts} give the two rank updates.

Interchanging $A$ and $B$ gives the same identity with $\beta$ in
place of $\alpha$. Apply it to $K'$ and undo the twist. For
$K'=\tau_AK$ it gives $y=x+(x+y)-2\rk\beta_{K'}$, so $\beta_{K'}$
surjects onto the rank-$x$ group $\HF(A,K')$. Thus
$R_{u_B}=0$ implies $\alpha_{K'}=0$ by \eqref{eq:composites}.
For $K'=\tau_B^{-1}K$ it gives $x=(x+y)+y-2\rk\beta_{K'}$, so
$\beta_{K'}$ is injective on the rank-$y$ group $\HF(B,K')$.
Now $R_{u_A}=0$ implies $\alpha_{K'}=0$. In both cases the relevant
top class acts on a group of rank $x+y$, so the last assertion follows
from Lemma~\ref{lem:top}.
\end{proof}

\subsection{The rank recurrence and the entropy bound}\label{sec:entropy}

\begin{lemma}\label{lem:seed}
For $K_r=\tau_{A_a}^r(B_b)$ with $r\ge1$, all edge maps
$\alpha_{ij,K_r}$ vanish. The Floer ranks are
$\dim\HF(A_i,K_r)=C_{ib}$ and $\dim\HF(B_j,K_r)=rC_{aj}$.
\end{lemma}
\begin{proof}
By \eqref{eq:twist-shifts}, $\HF(A_a,K_r)$ has rank one in degree
$r+1$. For every core $F$ disjoint from $A_a$, the twist sequence
identifies $\HF(F,K_r)$ with $\HF(F,B_b)$. This gives $x(K_r)$ and
leaves only the neighbors $B_j$ of $A_a$ to consider.

For $r=1$, the product
\[
 \HF^1(A_a,B_b)\otimes\HF^1(B_b,A_a)\longrightarrow\HF^2(B_b,B_b)
\]
is an isomorphism by duality. Thus $\HF(B_b,K_1)$ has rank one in
degree zero. For every other neighbor $B_j$, $\HF(B_j,B_b)=0$, so
the degree-two tensor term gives a single class in $\HF^1(B_j,K_1)$.
Suppose $\HF(B_j,K_r)$ has rank $r$ and vanishes in degrees above
$r$. The tensor term for the next twist is one-dimensional in
degree $r+2$. Its product therefore vanishes, and \eqref{eq:Seidel}
adds one class in degree $r+1$. Thus $\HF(B_b,K_r)$ has one class in
degree zero and one in each degree from $2$ to $r$. At any other
neighbor of $A_a$ the degrees are $1,\ldots,r$, each of rank one.
This proves the formula for $y(K_r)$.
It also gives $\alpha_{aj,K_r}=0$, since this map takes the
degree-$(r+1)$ generator to $\HF^{r+2}(B_j,K_r)=0$.

For $i\ne a$, let $\eta_F$ denote the map on $\HF(F,-)$ induced by
$B_b\to K_1\to\cdots\to K_r$. Since $A_i$ is disjoint from $A_a$,
$\eta_{A_i}$ is an isomorphism. Naturality gives
\[
 \alpha_{ij,K_r}\eta_{A_i}=\eta_{B_j}\alpha_{ij,B_b}.
\]
The right-hand side vanishes for $j\ne b$ because $\HF(B_j,B_b)=0$.
For $j=b$ the image of $\alpha_{ib,B_b}$ lies in $\Lambda u_{B_b}$,
which the map induced by $B_b\to K_1$ kills. This proves the
remaining vanishings.
\end{proof}

We now give the proof of Proposition~\ref{prop:rank-recurrence}.

\begin{proof}
At $K_*=K_4$, every edge map vanishes and $y_b(K_*)=4$ by
Lemma~\ref{lem:seed}. Suppose every edge map vanishes at $K$ and
$y_b(K)\ge4$. After the next twist, the new sphere $K'$ has
$y_b(K')\ge y_b(K)$ by Lemma~\ref{lem:edge} and \eqref{eq:Seidel}.
For every core $S$ we have
$\dim\HF(B_b,S)\le2$, so $K'$ cannot be isomorphic to a shift of
any core. Lemma~\ref{lem:top} now gives $R_{u_S}=0$ on $\HF(S,K')$
for every core $S$. Lemma~\ref{lem:edge} therefore gives
$\alpha_{ij,K'}=0$ on all edges incident to the twisting sphere.

For any other edge $A_i$--$B_j$, consider first a positive
$A$-twist. The map $K\to K'$ induces an isomorphism on $\HF(A_i,-)$,
since $A_i$ is disjoint from the twisting sphere. It commutes with
precomposition by the edge generator, so $\alpha_{ij,K}=0$ implies
$\alpha_{ij,K'}=0$. For a negative $B$-twist, use the map $K'\to K$.
This time it induces an isomorphism on $\HF(B_j,-)$, giving the same
implication. This completes the induction.

The $B$-twists thus replace $x$ by $x'=x+2qCy$ and preserve $y$.
The $A$-twists replace $y$ by $y+2pC^t x'$ and preserve $x'$.
Thus each application of $\Phi_{p,q}$ multiplies the Floer rank column
by the matrix in \eqref{eq:K3-matrix}. Iterating $n$ times from $K_*$
gives the formula.
\end{proof}

\begin{proof}[Proof of Theorem~\ref{main:k3}]
Let $M$ be the rank matrix in \eqref{eq:K3-matrix}. Reversing the
two invertible factors in
\[
 M=\begin{pmatrix}I&0\\2pC^t&I\end{pmatrix}
   \begin{pmatrix}I&2qC\\0&I\end{pmatrix}
\]
gives the two-block Penner matrix of Corollary~\ref{cor:penner-tree}
with exponents $2q$ and $-2p$. The two products are conjugate, so
\eqref{eq:penner-two-block} gives $\rho(M)=\lambda_\Gamma(p,q)$.

The matrix $M$ is primitive. Its diagonal is positive, and it has
positive entries in both directions along every edge of the connected
diagram $\Gamma$. By Proposition~\ref{prop:rank-recurrence}, the rank
column for $\Phi_{p,q}^nK_*$ is $M^n$ times the nonzero, nonnegative
initial column. The Perron--Frobenius theorem therefore gives
\[
 \lim_{n\to\infty}\lambda_\Gamma(p,q)^{-n}
       \dim_\Lambda\HF_X(L_v,\Phi_{p,q}^nK_*)>0
 \qquad\text{for every }v .
\]
Taking logarithms gives \eqref{eq:K3-growth}.

For $k>0$, apply Proposition~\ref{prop:floer-entropy} on $X$ with
the core spheres and $K_*$. Every $f$ in the class of $\Phi_{p,q}^k$
satisfies
\[
 \htop(f)\ge\lim_{n\to\infty}\frac1n\log\!\left(1+
       \sum_v\dim_\Lambda\HF_X(L_v,\Phi_{p,q}^{kn}K_*)\right)
       = k\log\lambda_\Gamma(p,q).
\]
For $k<0$ apply this estimate to $f^{-1}$ and use
$\htop(f^{-1})=\htop(f)$.

Seidel showed that the Dehn twist about the zero section of $D^*S^2$
generates its relative symplectic mapping class group $\Z$ and has
order two in the relative smooth mapping class group
\cite[Theorem~1 and Corollary~2]{SeidelCotangent}. The smooth isotopy
of its square to the identity fixes a boundary collar, so it extends
by the identity outside a cotangent neighborhood of each core sphere.
All exponents in $\Phi_{p,q}$ are even, so every power is smoothly
isotopic to the identity and acts trivially on $H^*(X;\Z)$.
\end{proof}

\begin{remark}
A holomorphic automorphism of a complex K3 surface acting trivially
on integral second cohomology is the identity
\cite[Proposition~15.2.1]{Huybrechts}.
For $k\ne0$, a holomorphic representative of $[\Phi_{p,q}^k]$ for
any complex K3 structure on $X$ would therefore be the identity and
have zero topological entropy.

An ADE configuration of smooth rational curves on a complex K3 gives
Lagrangian spheres for $\operatorname{Re}\Omega$, where $\Omega$ is
a nowhere-vanishing holomorphic two-form. Indeed, $\Omega$ vanishes
on each curve and $(\operatorname{Re}\Omega)^2=\Omega\wedge\bar\Omega/2$.
After hyperk\"ahler rotation, Theorem~\ref{main:k3} applies if the
rotated K\"ahler class is a positive multiple of a polarization
as in Definition~\ref{def:polarized-k3}.
\end{remark}

\subsection{The exact $A_2$ formula}\label{sec:recurrence}

Applying the rank recurrence from the proof of Theorem~\ref{main:k3}
to the first iterate of $L_2$ in type $A_2$ gives Corollary~\ref{cor:k3-A2}.
\begin{proof}[Proof of Corollary~\ref{cor:k3-A2}]
Let $A=L_1$, $B=L_2$, $r=2p$, $s=2q$, and $H_n=\Psi_{p,q}^{\,n}B$.
Write $x_n=\dim\HF(A,H_n)$ and $y_n=\dim\HF(B,H_n)$.
Since $\tau_B^{-s}B=B[s]$, the first iterate is $H_1=(\tau_A^rB)[s]$.
Lemma~\ref{lem:seed} gives $(x_1,y_1)=(1,r)$ and $\alpha_{H_1}=0$.

At $H_1$ the sum of the two ranks is $1+r\ge3$. Lemma~\ref{lem:edge}
therefore preserves the vanishing of $\alpha$ after the next twist,
and neither rank decreases. Iterating through the two twist blocks gives
\begin{equation}\label{eq:matrix}
 \binom{x_{n+1}}{y_{n+1}}
 =\begin{pmatrix}1&s\\r&1+rs\end{pmatrix}\binom{x_n}{y_n}
 \qquad(n\ge1).
\end{equation}
This matrix has determinant $1$ and trace $2+rs$. By
Cayley--Hamilton,
\begin{equation}\label{eq:scalarrec}
 x_{n+2}=(2+rs)x_{n+1}-x_n\qquad(n\ge0),
\end{equation}
where the case $n=0$ follows from $x_0=x_1=1$ and $x_2=1+rs$.
The characteristic roots are $\lambda_{p,q}^{\pm1}$. The expression
$D_n(p,q)$ in \eqref{eq:lambda} satisfies this recurrence and both
initial conditions, which proves \eqref{eq:growth}.
\end{proof}

\appendix
\section{Proof of the global twist sequence}
\label{app:seidel-comparison}

Proposition~\ref{prop:closed-sequence} follows from Seidel's twist
triangle for projective Calabi--Yau surfaces. We identify its evaluation
map with the Floer product in $X$ and its coefficients with those of
Section~\ref{sec:k3-main}.

Put $\Omega=\omega/c$, so that $[\Omega]=c_1(\mathcal H)$, and choose
a unitary connection on $\mathcal H$ of curvature $-2\pi i\Omega$.
A Lagrangian is \emph{rational} if the induced flat connection has
finite holonomy \cite[\S\,6a]{SeidelQuartic}. On a Lagrangian sphere
the holonomy is trivial. We may therefore choose a parallel unit
section for the flat multisection in Seidel's brane data. His
coefficient field is $\mathbb K=\bigcup_{m\ge1}\C((q^{1/m}))$.
Write $\Hom^k(P,Q)=\Hom_{D^b\mathcal F(X)}(P,Q[k])$.

\begin{theorem}[{\cite[\S\,8c, Lemmas~8.4 and~8.14]{SeidelQuartic}}]
\label{app:seidel-closed-fukaya}
For each graded Lagrangian sphere $S\subset X$ there is a compatible
almost complex structure $J$ for which $X$ contains no nonconstant
$J$-holomorphic spheres and $S$ bounds no nonconstant $J$-holomorphic
discs. With these boundary data $S$ defines an object of $\mathcal F(X)$
up to quasi-isomorphism. Graded isotopic spheres give quasi-isomorphic
objects. Their morphism spaces are Floer cohomology groups, and
composition is the Floer triangle product in $X$.
\end{theorem}

\begin{theorem}[{
\cite[(2.14)--(2.15)]{SeidelQuartic}}]
\label{app:seidel-algebraic-twist}
If $\Hom^*(P,Q)$ has finite total dimension, the algebraic twist
$T_P(Q)$ fits into the distinguished triangle
\[
 \Hom^*(P,Q)\otimes P\xrightarrow{\mathrm{ev}}Q
 \longrightarrow T_P(Q)\longrightarrow(\Hom^*(P,Q)\otimes P)[1].
\]
In a homogeneous basis $(y_j)$ of $\Hom^*(P,Q)$, the component of
$\mathrm{ev}$ on $P[-|y_j|]$ represents $y_j$.
\end{theorem}

\begin{theorem}[{
\cite[\S\,9c, Proposition~9.4]{SeidelQuartic}}]
\label{app:seidel-geometric-twist}
Let $L$ be a graded Lagrangian sphere, $Q$ a rational brane, and
$\widetilde{\tau}_L$ the graded Dehn twist normalized to be the
identity outside its support. Then
$T_L(\widetilde{\tau}_L^{-1}Q)\cong Q$ in $D^b\mathcal F(X)$,
where the inverse image brane carries the pulled-back boundary data.
\end{theorem}

\begin{proof}[Proof of Proposition~\ref{prop:closed-sequence}]
The same graded twist represents $\tau_L$ for $\omega$ and for
$\Omega$. Indeed, if $\psi$ is a cotangent chart for $\omega$, then
$(u,p)\mapsto\psi(u,cp)$ is a chart for $\Omega$, with twist profile
$r\mapsto\vartheta(cr)$ in place of $\vartheta(r)$. The twist also
lifts to a connection-preserving map on $\mathcal H$. The induced
connection on $\Hom(\mathcal H,\tau_L^*\mathcal H)$ is flat, and
$X$ is simply connected.

Take $Q=\widetilde{\tau}_L K$ in
Theorem~\ref{app:seidel-geometric-twist}. The pulled-back sphere
represents $K$ by Theorem~\ref{app:seidel-closed-fukaya}, so
$T_L(K)\cong\tau_LK$ with the induced grading. Graded Hamiltonian
pushoffs make the spheres transverse, giving finite-dimensional
Floer cochains. Applying $\Hom^*(F,-)$ to the triangle of
Theorem~\ref{app:seidel-algebraic-twist} gives the long exact sequence
over $\mathbb K$. On the summand indexed by $y_j$, evaluation sends
$x$ to $y_j\circ x$. Thus its map on the tensor product is the
Floer product in $X$.

It remains to identify the coefficients. For each Floer generator $x$,
choose a real lift $a_x$ of its action modulo $\Q$. Its orientation
line has coefficients in $q^{a_x}\mathbb K$. A rigid strip or triangle
with inputs $x_i$ and output $x_0$ has weight $\pm q^{E}$, where
$E+\sum_i a_{x_i}-a_{x_0}\in\Q$. Here $E$ is the Floer energy for
a strip and the energy minus the integrated curvature of the
perturbation for a triangle \cite[(8.18)--(8.20)]{SeidelQuartic}.
For unperturbed data $E$ is the $\Omega$-area.

Extend scalars by $q^r\mapsto T^{cr}$ and identify the generator
$q^{a_x}$ of each coefficient line with $T^{ca_x}x$. In the
intersection basis the weight becomes
\[
 \pm T^{c(E+\sum_i a_{x_i}-a_{x_0})}
       T^{ca_{x_0}-c\sum_i a_{x_i}}=\pm T^{cE}.
\]
Scaling $\Omega$ and the Hamiltonians by $c$ preserves the perturbed
curves and scales their energies and curvature terms by $c$.
Thus the differential and product have the $\omega$-weights.
Extension to the field $\Lambda$ is exact, giving \eqref{eq:Seidel}.

The second arrow of the triangle gives $\iota:K\to\tau_LK$.
Applying the same triangle to $\tau_L^{-1}K$ gives
$\pi:\tau_L^{-1}K\to K$. Postcomposition with either arrow commutes
with precomposition by associativity, as required in
Proposition~\ref{prop:rank-recurrence}.
\end{proof}

\section{Faithfulness of slice monodromy in types \texorpdfstring{$D$ and $E$}{D and E}}
\label{app:slice-faithful-de}

By Lemma~\ref{lem:slice-monodromy}, each simple meridian acts by a Dehn
twist about a vanishing sphere. To distinguish the resulting words, we
apply the theorem of Brav and Thomas to the graded vanishing spheres.

\begin{theorem}[{Seidel \cite[Corollary~17.17]{SeidelBook}}]
\label{thm:appendix-exact-twist}
Let $L$ be a framed exact graded Lagrangian $2$-sphere with spin structure
in a graded Liouville four-manifold. For every compact exact graded brane $Q$,
the normalized graded Dehn twist satisfies
\[
 \widetilde{\tau}_L Q\cong T_L(Q)
 =\operatorname{Cone}\bigl(
 \operatorname{RHom}(L,Q)\otimes L\xrightarrow{\mathrm{ev}}Q\bigr)
\]
in the derived Fukaya category.
\end{theorem}

\begin{theorem}[{Brav--Thomas \cite[Proposition~3.1, Theorem~3.1 and its proof]{BravThomas}}]
\label{thm:appendix-brav-thomas}
Let $(S_i)$ be a configuration of $2$-spherical objects with ADE diagram
$\Gamma$ in an enhanced linear triangulated category. Their spherical twists give a
faithful action of $\Art_\Gamma$. For $S=\bigoplus_i S_i$, a positive
Artin element $\alpha$ is determined by the isomorphism class of
$T_\alpha S$.
\end{theorem}

\begin{lemma}[{Seidel \cite[\S\,2 and Lemmas~5.6--5.7]{SeidelGraded}}]
\label{thm:appendix-grading-lift}
The normalized lift of a Dehn twist about a graded Lagrangian
$2$-sphere $L$ is the identity outside the twist neighborhood and
sends $L$ to $L[-1]$. On a connected manifold, two graded lifts of
the same symplectomorphism that agree on a boundary collar coincide.
\end{lemma}

\begin{lemma}\label{lem:slice-faithful}
For $\Gamma$ of type $D$ or $E$, the slice monodromy
$\rho_{\mathrm{slice}}\colon\Art_\Gamma\to\pi_0\Symp_c(\Int F_\Gamma)$
is injective.
\end{lemma}
\begin{proof}
Use the grading defined by the residue holomorphic two-form on
$F_\Gamma$. Work over $\C$ in the full triangulated subcategory of the
derived exact Fukaya category split-generated by the vanishing spheres
$L_i$, each with its unique spin structure. Root the tree and shift the
gradings successively so that each intersection from parent to child
has degree one. The reverse intersection also has degree one by
duality. Thus $\HF^*(L_i,L_j)$ is one-dimensional in degree one at
each edge and zero for distinct nonadjacent vertices. Exactness gives
$\HF^*(L_i,L_i)\cong H^*(S^2;\C)$. Floer composition and trace give
the perfect degree-two pairing on compact exact branes
\cite[Lemma~15]{GanatraCY}, which extends to shifts, cones and direct
summands. Hence the $L_i$ form a $\Gamma$-configuration of
$2$-spherical objects.

Suppose $\rho_{\mathrm{slice}}(\beta)=1$. The geometric twist word is
then symplectically isotopic to the identity through maps fixing a
common boundary collar. The fiber retracts onto its tree of spheres,
so $H_c^1(\Int F_\Gamma;\R)\cong H_3(F_\Gamma;\R)=0$ and the isotopy
is Hamiltonian. Lift it starting at the graded identity. The lift
fixes the collar throughout, so its endpoint is the product of the
normalized graded twists by Lemma~\ref{thm:appendix-grading-lift}.
Graded Hamiltonian continuation therefore identifies the image of
each $L_i$ with $L_i$ itself.

Put $S=\bigoplus_i L_i$. Theorem~\ref{thm:appendix-exact-twist}
identifies each graded geometric twist with the spherical twist on
objects. Applying the same formula to $\widetilde{\tau}_{L_i}^{-1}Q$
gives the identification for inverse twists. Thus $T_\beta S\cong S$.
Write $\beta=\alpha^{-1}\gamma$ with $\alpha,\gamma$ positive
\cite[\S\,2]{BravThomas}. Then $T_\gamma S\cong T_\alpha S$.
Theorem~\ref{thm:appendix-brav-thomas} gives $\alpha=\gamma$, so
$\beta=1$.
\end{proof}

\section{Arnold's questions addressed here}
\label{subsec:arnold-problems}

We use the problem numbers and page references of \cite{ArnoldProblems}.

\subsection*{Flat transport and braid relations}
Problem~1991-12 (p.~80) asks whether the monodromy representation of an
$A_n$ plane-curve singularity lifts to the diffeomorphism group of the fiber, and asks the
same question in higher dimensions and for symplectomorphisms.
Problem~1993-13 (p.~87) asks for representatives $a,b$ of the meridian
and longitude Dehn twists on a torus with one boundary component such
that $aba=bab$ as maps, preferably fixing the boundary pointwise.

For the transport problem of \cite[p.~101]{ArnoldMilnor}, Shmelev
constructed flat symplectic transport with Weyl group holonomy on the
classical matrix families of types $A$, $B$, $C$ and $D$ with distinct
eigenvalues. He also proved that in rank greater than one no connection
on these families can have trivial holonomy on the sphere bundle at
infinity.\footnote{A.~S.~Shmelev,
\emph{Symplectic connections on the fiberings of regular orbits under
the adjoint representations of simple Lie groups}, Funct. Anal. Appl.
\textbf{30} (1996), 36--41, \S\,2 and Theorem~3.5,
\href{https://doi.org/10.1007/BF02509555}{doi:10.1007/BF02509555}.}
For any connected complex semisimple group, Theorem~\ref{main:flat-orbits}
gives a complete flat connection with Weyl group holonomy and the
prescribed form on a reference orbit with real adjoint spectrum.
For the flag and slice actions,
Theorems~\ref{main:flag} and~\ref{main:slice} give the braid relation
$[a][b][a]=[b][a][b]$ for twists at adjacent vertices in the
respective symplectic mapping class groups.

Problem~1994-9 (p.~99) asks whether the $A_2$ surface family
\[
 \{(x,y,z)\in\C^3:x^3+ux+v+y^2+z^2=0\},
 \qquad 4u^3+27v^2\ne0,
\]
admits a flat symplectic connection. Theorem~\ref{main:slice} determines
its monodromy on compact Milnor fibers up to symplectic isotopy.
Restricting the flat connection of Theorem~\ref{main:flat-orbits} would
require its horizontal spaces to be tangent to the slice. The existence
of a flat symplectic connection on this family remains a separate question.

\subsection*{The full flag and its homology}
\begin{arnoldproblem}[1992-2, p.~81, flag monodromy]
Study the braid action arising from the coadjoint representation of
$SL_N(\C)$ on the full flag manifold and its tangent bundle,
identified as a real vector bundle with the cotangent bundle by a
Riemannian metric. Construct a theory of monodromy and
variation for non-isolated singularities, with boundary conditions
along the strata that meet the boundary of a Milnor ball. The theory
should apply to the characteristic-polynomial map
$\chi_N(x)=\det(zI-x)$ and extend the Brieskorn--Grothendieck
description of transverse ADE surface singularities.
\end{arnoldproblem}
The trace pairing identifies the coadjoint and adjoint representations.
Theorem~\ref{main:flag} gives the resulting action on $T^*\Fl_N$ with
the prescribed behavior at infinity, and Theorem~\ref{main:slice} treats
the subregular transverse slices. A theory of variation along all the
singular strata would also have to specify how relative cycles change
under transport near each stratum.

\begin{arnoldproblem}[1993-12, p.~87, the homology action]
Describe the braid action on the integral homology of a regular orbit
$X_b\cong T^*\Fl_N$, including its restriction to the subgroups
associated with non-isolated singularities of the fibers.
\end{arnoldproblem}
\begin{remark}\label{rmk:forget-infinity}
In the full symplectomorphism group, the profiles
$(1-s)\vartheta+s\pi/2$, $0\le s\le1$, deform $\tau_i$ to the Weyl
cotangent lift \eqref{eq:zero-extension}. This lift has square $\id$.
In the $SU(N)$-equivariant group with the prescribed condition at
infinity, Jin's braid homomorphism sends $\tau_i^2$ to $\sigma_i^2\ne1$
\cite[Theorem~1.4 and \S\,3.2]{Jin}.
\end{remark}

\phantomsection\label{subsec:homological-action}
By Theorem~\ref{main:flat-orbits} and Remark~\ref{rmk:forget-infinity},
both actions on homology factor through the permutation map
$Br_N\to S_N$. On the zero section, $\tau_i$ exchanges the $i$th and
$(i+1)$st orthogonal lines. Thus for the tautological classes
$x_j=c_1(V_j/V_{j-1})$ and the transposition $s_i=(i\ \ i+1)$,
\[
                    (\tau_i^{-1})^*x_j=x_{s_i(j)} .
\]
Every subgroup acts through its permutation image. In particular, pure
braids act trivially on integral homology and cohomology.

\subsection*{Floer ranks and iterated intersections}
Problem~1974-6 (p.~20) asks which topological theories have symplectic
counterparts. Theorem~\ref{main:k3} gives an example on a polarized K3
surface with an $A_2$ configuration. For every nonzero integer $m$,
the twist word $\Psi_{p,q}^{\,m}$ is smoothly isotopic to the identity.
Every representative $f$ of its symplectic mapping class satisfies
\[
                 \htop(f)\ge |m|\log\lambda_{p,q}>0 .
\]

\begin{arnoldproblem}[1994-48, p.~107, iterated intersections]
Let $X_0,Y_0$ be compact submanifolds of a compact manifold $M$, and
let $f\colon M\to M$ be differentiable. How does the topological
complexity of $f^n(X_0)\cap Y_0$ grow with $n$? For smooth data,
is it almost always bounded above by an exponential in $n$? One
measure Arnold proposes is the sum of the Betti numbers of the
intersection.
\end{arnoldproblem}
For the spheres in Corollary~\ref{cor:k3-A2}, this sum counts points
whenever the intersection is transverse. If $f=\Psi_{p,q}$ and
$\eta$ is any Hamiltonian diffeomorphism such that $\eta(f^n L_2)$
meets $L_1$ transversely, then
\[
 \#\bigl(\eta(f^n L_2)\cap L_1\bigr)
 \ge\dim_\Lambda\HF_X(L_1,f^n L_2)=D_n(p,q).
\]
By \eqref{eq:lambda}, this lower bound is asymptotic to
$\lambda_{p,q}^n/(\lambda_{p,q}+1)$. For an upper bound on the Floer ranks of iterated
Lagrangian intersections, see Seidel's account of Arnold's estimate
\cite[Proposition~5.11]{SeidelLecture}.

\end{document}